\documentclass[11pt]{amsart}
\usepackage{amsmath,amssymb} 
\usepackage{latexsym}
\newtheorem{thm}{Theorem}[section]
\newtheorem{cor}[thm]{Corollary}

\newtheorem{lem}[thm]{Lemma}
\newtheorem{prop}[thm]{Proposition}
\theoremstyle{definition}

\theoremstyle{remark}
\newtheorem{rem}{Remark}[section]
\newtheorem{ex}{Example}[section]
\numberwithin{equation}{section}
\newcommand{\thmref}[1]{Theorem~\ref{#1}}
\newcommand{\secref}[1]{\S\ref{#1}}
\newcommand{\lemref}[1]{Lemma~\ref{#1}}
\newcommand{\propref}[1]{Proposition~\ref{#1}}
\newcommand{\corref}[1]{Corollary~\ref{#1}}
\newcommand{\remref}[1]{Remark~\ref{#1}}
\newcommand{\disp}{\displaystyle}
\newcommand{\nc}{\newcommand}
\nc{\bib}{\bibitem}
\nc{\on}{\operatorname}
\nc{\res}{\operatornamewithlimits{Res}}
\nc{\ahat}{\on{Ahat}}
\nc{\ul}{2\pi\sqrt{-1}\,}
\nc{\arr}{\rightarrow}
\nc{\al}{\alpha}
\nc{\C}{{\mathbb C}}
\nc{\Cn}{{\mathbb C}^n}
\nc{\Z}{{\mathbb Z}}
\nc{\R}{{\mathbb R}}
\nc{\bB}{\mathbf{B}}
\nc{\bb}{\mathbf{b}}
\nc{\A}{{\mathfrak A}}
\nc{\B}{{\mathfrak B}}
\nc{\mC}{{\mathfrak C}}
\nc{\fA}{\widehat\A}
\nc{\I}{{\mathfrak I}}
\nc{\eps}{\varepsilon}
\nc{\M}{{\mathfrak M}}
\nc{\ml}{\ll}
\nc{\mgg}{\M^G_g}
\nc{\vs}{V^*}
\nc{\Sym}{\on{Sym}}
\nc{\kh}{\eta}

\renewcommand{\H}{\on{H}}
\nc{\Hom}{\mathrm{Hom}}
\nc{\tensor}{\otimes}
\nc{\spf}{$\operatorname{Spin}(5)$~}
\nc{\suk}{$\operatorname{SU}(2)$~}
\nc{\suh}{$\operatorname{SU}(3)$~}
\nc{\sun}{$\operatorname{SU}(n)$~}
\nc{\ssum}{\sum}
\nc{\im}{\on{im}}
\nc{\rat}{R_\A}
\nc{\merom}{M_\A}
\nc{\flag}{\mathrm{Flag}}
\nc{\bcb}{\on{NBC}}
\nc{\nbc}{{\bf nbc}}
\nc{\End}{\on{End}}
\nc{\id}{\mathrm{id}}
\nc{\ach}{\A^*}  
\nc{\starr}{$*$}
\nc{\Or}{\A^n_\mathrm{ind}} 
\nc{\Rf}{\Z \Or} 
\nc{\OS}{\mathbf{OS}}
\nc{\hm}{\hat M}
\nc{\ho}{\hat O}
\nc{\ha}{\hat{\mathfrak A}}
\nc{\hos}{\hat \mu}
\nc{\vex}{\bx}
\nc{\vey}{\by}
\nc{\vol}{\on{vol}}
\nc{\cube}{\Box}
\nc{\oma}{\Omega_{\A}}
\nc{\less}{\backslash}
\nc{\mer}{M_\A}
\nc{\Cone}{\on{Cone}}
\nc{\cla}{\mathbb{E}(\Lambda)}
\nc{\bx}{\mathbf{x}}
\nc{\bbx}{\mathbf{X}}
\nc{\by}{\widehat{\mathbf{y}}}
\nc{\bq}{\widehat{\ba}}
\nc{\rt}{\hat R_\A}
\nc{\ttt}{\tilde t}
\nc{\dn}{{\mathfrak D}_n}
\nc{\bbb}{BB^n}
\nc{\bara}{\overline{|\A|}}
\nc{\abez}{\A\backslash H}
\nc{\diag}{\mathrm{diag}}
\nc{\diagd}{\widehat{\mathrm{diag}}}
\nc{\bbeta}{\mathbf{a}}
\nc{\mudel}{\hat{\mu}_\Delta}
\nc{\qdel}{\hat{q}_\Delta}
\nc{\hatir}{\widehat{\ir}_\Delta}
\nc{\bz}{\vec{\mathbf{z}}}
\nc{\ires}{\operatornamewithlimits{IRes}}
\nc{\Gs}{\Gamma^*}
\nc{\Gc}{\check\Gamma}
\nc{\Tc}{\check\Theta}
\nc{\Ts}{\Theta^*}
\nc{\Ch}{\C_h}
\nc{\ver}{\mathrm{Ver}}
\nc{\bla}{\vec{\lambda}}
\nc{\todd}[1]{\frac{1}{1-e^{#1}}}
\nc{\fel}{\frac12}
\nc{\bfi}{\mathbf{\phi}}
\nc{\ac}{\A}
\nc{\compl}{U(\A)}
\nc{\ct}{\mathrm{CT}}
\nc{\cta}{\ct^\A}
\nc{\ctp}{\underset{p}{\mathrm{CT}}}
\nc{\perm}{\mathcal{S}}
\nc{\aind}{\A^n_\mathrm{ind}}
\nc{\td}{\mathrm{Todd}}
\nc{\ict}[1]{\operatornamewithlimits{iCT}_{#1}}
\nc{\tq}{\tilde q^\Gamma_\tau}
\nc{\qh}{{\hat Q}}
\nc{\vr}{V_\R}
\nc{\vdr}{V^*_\R}
\nc{\PP}{\mathrm{Pr}}
\nc{\lig}{\mathfrak{g}}
\nc{\lt}{\mathfrak{t}}
\nc{\tdual}{\lt^*}
\nc{\dplus}{|\Delta^+|}
\nc{\irr}{\mathrm{Irrep}}
\nc{\xh}{{\hat{x}}}
\nc{\frh}{\mathfrak{h}}
\nc{\isom}{\cong}
\nc{\orb}{\mathrm{Conj}}
\nc{\ar}{\mathfrak{R}}
\nc{\rfr}{R_\ar}
\nc{\trfr}{R_\ar[\Theta]}
\nc{\gz}{\mathfrak{z}}
\nc{\dc}{\mathfrak{C}}
\nc{\tom}{\tilde\omega}
\nc{\smn}{\sum_{i=0}^n}
\nc{\dch}{\mathfrak{C}^{\leq h}}
\nc{\fs}[1]{\overline{1,#1}}
\nc{\bet}{{e}_t}
\nc{\bey}{{e}_y}
\nc{\bew}{{e}_w}
\nc{\ctc}{\C[\check\Theta]}
\nc{\catc}{\C_\A[\check\Theta]}
\nc{\catcm}{\C_\A^\mu[\check\Theta]}
\nc{\ctatc}{\C_{\tilde \A}[\check\Theta]}
\nc{\ctatcm}{\C^\mu_{\tilde \A}[\check\Theta]}
\nc{\be}{{e}}
\nc{\amh}{\A\backslash H}
\nc{\arh}{\A_{|H}}
\nc{\tA}{\tilde\A}
\nc{\oset}{\subset\!\!\!\!\prec}
\nc{\um}{\underline m}
\nc{\un}{\underline n}
\nc{\ala}{\underline}
\nc{\ba}{\mathbf{a}}
\nc{\Ac}{\A^\circ}
\nc{\Ahc}{\widehat{\A}^\circ}
\renewcommand{\tt}[1]{\left(1-e^{#1}\right)}
\newcommand{\cx}{x^\circ}
\nc{\ty}{\tilde y}
\nc{\cy}{y^\circ}
\nc{\bD}{\bar D}
\nc{\fB}{\widehat{\B}}
\nc{\inter}{\mathrm{int}}
\nc{\BB}{P}
\nc{\bah}{\widehat{\ba}}
\nc{\ctt}{\widetilde{\mathrm{CT}}}
\renewcommand{\vert}{\mathrm{vx}}
\nc{\har}{\widehat{\ar}}
\nc{\sign}{\mathrm{sign}}
\nc{\At}{\A^{+\Theta}}
\nc{\Atp}{\A_p^{+\Theta}}
\nc{\LL}{\mathfrak{L}}

\begin{document}
\title[A residue theorem for rational
trigonometric sums] {A residue theorem for rational trigonometric
sums and Verlinde's formula}
\author{Andr\'as Szenes}
\thanks{The research was supported by NSF grant DMS-9870053 and NSA
  grant \#6800900}
\address{Massachusetts Institute of Technology, Department of Mathematics}
\email{szenes@math.mit.edu}

\maketitle

\section{Introduction}

The central objects of the present work are rational trigonometric
sums such as
\begin{equation} \label{eq:verex}
\sum\frac1%
{\sin^2\frac{\pi m}k\sin^2\frac{\pi n}k\sin^2\frac{\pi (m+n)}k},
\quad m,n\in\Z, \,0<m,n<k,\, m+n\neq k,
\end{equation}
where $k$ is a fixed positive integer.

The interest in such sums was motivated by a beautiful formula of
E.~Verlinde for the dimension of the ``space of conformal
blocks'' in the WZW theory \cite{Ver}.  The data for Verlinde's
formula is a simple simply-connected Lie group $G$, non-negative
integers $g$ and $k$, and, in the simplest case of one puncture, a
dominant highest weight $\lambda$ of $G$ satisfying certain
conditions. The result is a non-negative integer, which we denote by
$\ver_g(\lambda;k)$. The sum in \eqref{eq:verex} is  an example of
Verlinde's expression; up to some normalization, it represents the
case of $G=SU(3)$, $g=2$, $\lambda=0$.

This formula turned out to have a close relationship with the topology
of the moduli spaces of flat connections over Riemann
surfaces. Indeed, under certain assumptions, 
$\ver_g(k\lambda;k)$ is expected to agree with the Hilbert polynomial 
\[ \int_\M e^{kc_1(\LL)} \td(\M) \]
of a certain genus-$g$ moduli space of flat $G$-connections $\M$
endowed with a line bundle $\LL$, both depending on the data
$(G,g,\lambda)$. The discovery of this fantastic ``coincidence''
opened the way to computing various intersection numbers on the moduli
spaces \cite{thad}; such computations looked utterly impossible until
then \cite{AB}.

Looking at \eqref{eq:verex}, one might note that it is not at all
clear that the value of this finite sum is polynomial in $k$, as would
follow from the agreement with the Hilbert polynomial of a space. 
To uncover the topological information hidden in Verlinde's formula,
one needs to find a calculus that replaces the sum smeared in space by a
compact expression which is manifestly a polynomial in $k$.  A similar
problem arose with the evaluation of rational sums such as
\[ \sum \frac{1}{ m^2 n^2( m+ n)^2},
\quad m,n,m+n\in \Z^{\neq0}.
\]
These sums appeared in the work of Witten on 2-dimensional gauge
theory \cite{wit1,wit2}, and they again turned out to have a close
relationship  to the topology of the abovementioned moduli spaces. The
exact evaluation of these sums was left open in \cite{wit2}.

A solution to these problems was suggested by the author in
\cite{comb}.  We conjectured that there exist certain local
functionals on the space of rational and rational trigonometric sums
corresponding to the Lie group G, which would enable one to localize
both types of sums, and, moreover, that these functionals would
coincide in the two cases in such a way as to provide a Riemann-Roch
calculus on the moduli spaces. The functional was described in
\cite{comb} for the case of $G=SU(n)$ explicitly. Surprisingly, it had
the form of a single {\em iterated residue}.  This argument and the
functional, in a somewhat modified form, was later used by Jeffrey and
Kirwan \cite{JK} to give a computation of the Hilbert polynomial in
this case.

The functional for the case of $G=SU(n)$ is very simple, but in the
case of other groups a similarly symmetric formula does not seem to
exist. So we approached this problem from a more general point of
view -- from the point of view of arbitrary hyperplane arrangements, of
which the Stiefel diagram of a Lie group is a particular example. We
described the abovementioned functional in \cite{bern} for the case of
arbitrary rational sums. The present paper gives the answer for the rational
trigonometric case. 

Our effort was strongly motivated by the work of Bismut and Labourie
\cite{BL}. They computed the Hilbert polynomial of the moduli space
for an arbitrary group in terms of rational sums, but, curiously,
their formula did not seem to coincide with Verlinde's expression. One
of the main goals of this work was to prove the equality of the two
expressions. We achieved this in most cases, but there is still a
bit of mystery left when $g=0$. This will be duly explained in the
main body of the paper.

Finally, we need to mention a circle problems closely related to
computing rational trigonometric sums: the problems of partition
functions and counting lattice points in polytopes. Our localization
theorem is somewhat analogous to the results of Brion and Vergne on
vector partition function \cite{BVvp}. In fact, the exact relation is
worth investigating further \cite{SzV}.

The contents of the paper are as follows: in \S2 we recall the results
of \cite{bern} on rational sums and Bernoulli polynomials
corresponding to central hyperplane arrangements; in \S3 we extend
these results to affine hyperplane arrangements. The main theorem,
\thmref{thm:main-trig}, which gives a local formula for rational
trigonometric sums, is given in \S4 and the application to the formula
of Bismut and Labourie is detailed in \S5.

\textbf{Acknowledgments} The author is grateful to Jean-Michel Bismut
and Michel\`e Vergne, whose support and warm hospitality made this
work possible. Noam Elkies, Pavel Etingof and Victor Kac supplied some
crucial advice.

\section{Central hyperplane arrangements}
\label{sec:central}

In this section we review the results of \cite{bern} (cf. \cite{BV}
for an alternative treatment). 

\subsection{Notation and Conventions}
Let $\A$ be a central and essential hyperplane arrangement (HPA) in an
$n$-dimensional complex vector space $V$, i.e.  a collection of
hyperplanes in $V$ such that $\cap\A=\{0\}$.  Denote by $\rat$ the
rational functions on $V$ with poles along $\cup\A$, by $\mer$ the
meromorphic functions defined in a neighborhood of 0 with poles along
$\cup\A$, and by $\compl$ the complement $V\setminus\cup\A$.

To simplify our notation, we impose a linear ordering relation $\prec$
on $\A$ and assume that the hyperplanes are indexed accordingly:
$(\A,\prec)=(H_1,\dots,H_N)$. We will make this ordering explicit in
the notation whenever it is used in our constructions in an essential
manner. If an $m$-element subset $\ba$ of $\A$ is ordered, we will
write $\ba\in\A^m$, and think of it as of a sequence of elements of
$\A$. If the ordering of $\ba$ is consistent with $\prec$, then we
will write $\ba\oset\A$.  Often it will be convenient to choose a
linear form $x_i$ for each hyperplane $H_i\in\A$ such that
$H_i=\{x_i=0\}$.  We will use the notation $\fA$ for this ordered
set of linear forms.  Again, we will use the corresponding notations
$\bq\in\fA^m$ and $\bq\oset\fA$. Finally, we will denote the $i$th
element of $\ba$ by $H_{i,\ba}$, the $i$th element of $\bq$ by
$x_{i,\bq}$ and the function $e^{2\pi \sqrt{-1}x_{i,\bq}}$ by
$e_{i,\bq}$.

We will call a set $\{H_i\}_{i=1}^m$ of $m$ hyperplanes in $V$
{\em independent} if $\dim \cap H_i = n-m$. This is equivalent to
saying that the corresponding linear forms are linearly
independent.

An important convention throughout the paper is that underlining a
symbol means multiplying it by $2\pi\sqrt{-1}$. For $y\in V^*$ we will
write $e_y$ for the function $e^{2\pi\sqrt{-1}\,y}$ on $V$, which thus
may also be written as $e^{\ala y}$.

Finally, the notation $\fs k$ will sometimes be used to denote the
set of first $k$ natural numbers.

\subsection{The Constant Term}

One can associate to every hyperplane $H_i\in\A$ a closed
holomorphic differential 1-form $\alpha_i=dx_i/x_i$ on $\compl$,
where $\{x_i=0\}=H_i$. These 1-forms are called {\em logarithmic
  differential forms}; as they are homogeneous in the linear form
$x_i$, they do not depend on its choice. The key fact of the
topology of $\compl$ is the existence of an injective ring
homomorphism $q:H^*(\compl,\C) \rightarrow \Omega^*(\compl)$,
which assigns to every cohomology class a closed holomorphic
differential form corresponding to it in de Rham theory.  The
image $\Omega^*_\A=\im(q)$ of the cohomology ring is generated in
degree 1 by the logarithmic 1-forms $\{\alpha_i,\,
i=1,\dots,N\}$.

The map $q$ allows us to define a generalized \emph{constant
  term} functional $\cta:\mer\rightarrow \C$, which is of degree
0 with respect to the natural grading on $\rat\subset\mer$. Given
a representation $\sum_i Z_i\tensor \beta^i\in
H_n(\compl,\C)\tensor H^n(\compl,\C)$ of the invariantly defined
diagonal element derived from the natural complete pairing of
$H_n(\compl,\C)$ and $H^n(\compl,\C)$, we may form the functional
\begin{equation}
  \label{eq:cta}
\ct^\A : f\mapsto \sum_i \int_{Z_i}f\,q(\beta^i).   
\end{equation}
This functional is invariantly defined; it depends solely on the
HPA $\A$, with no additional choices made. 

If $|\A|=n=\dim V$, i.e. the HPA is simple, then this functional is
equal to the constant term of the Laurent expansion of the function
$f$ near 0.  In this case we will simply write $\ct$ instead of
$\cta$.  When $|\A|>n$, a more involved algebraic computational device
is available; this will be detailed below.

\subsection{Iterated constant term functionals}  Let
$\aind\subset\A^n$ be the set of independent ordered $n$-tuples of
hyperplanes in $\A$.  Given $\ba\in\aind$ and a permutation
$\tau\in\perm_n$, denote by $\ba^\tau$ the element
$(H_{\tau(1),\ba},\dots,H_{\tau(n),\ba})\in\aind$.  Then it follows
from the description and properties of the map $q$ that the linear
space $\Omega^n_\A=q(H^n(\compl,\C))$ is spanned by the forms
$$\alpha_\ba=\alpha_{1,\ba}\wedge\dots\wedge\alpha_{n,\ba}$$ 
as $\ba$ varies in $\aind$.

Every $\ba\in\aind$ defines an \emph{iterated constant term} functional
\[ \ict
\ba=\ct_{H_{1,\ba}}\ct_{H_{2,\ba}}\dots\ct_{H_{n,\ba}}:
M_\A\longrightarrow\C,\] which is obtained by sequentially
applying the 1-dimensional constant term functional with respect
to each of the hyperplanes in $\ba$, while keeping the preceding
variables constant.  More precisely, the symbol
$\ct_{H_{n,\ba}}f$ means taking the 1-dimensional constant term
of $f$ along each of the lines $\{x_{i,\ba}=a_i,\, i=1,2,\dots
n-1\}$, where $(a_1,\dots,a_{n-1})$ is a fixed $(n-1)$-tuple of
complex numbers. Considering $(a_1,\dots,a_{n-1})$ to be
coordinates on $H_{n,\ba}$, we can think of $\ct_{H_{n,\ba}}f$ as
of a function on $H_{n,\ba}$. Then we replace $V$ with
$H_{n,\ba}$ and the function $f$ with the function
$\ct_{H_{n,\ba}}f$ on $H_{n,\ba}$, and we continue the process,
taking $\ct_{H_{n-1,\ba}}$ of $\ct_{H_{n,\ba}}f$, etc, finally
arriving at a number. 

We should point out that the notation $\ct_{H_{n,\ba}}f$ is
somewhat misleading since this quantity depends on the rest of
the hyperplanes as well; this is similar to the problem with the
notation used for partial derivatives.  The iterated constant
term functional $\ict \ba$ usually depends on the order of the
hyperplanes in $\ba$ (cf. \cite{bern} for further details). To
simplify our notation, we will often write $\ct_{x}$ instead of
$\ct_H$ if the appropriate forms have been introduced.  Finally,
while the above definition of iterated constant terms might seem
somewhat complicated, computationally this procedure is very
simple, and is, in fact, a built-in function of software packages
such as \textsc{Maple}. It was through such computer
experimentation that the author became acquainted with the
concept.
\begin{ex} Let 
$$ f = \frac{e^{x-2y}}{xy(x+y)}.$$
Then
$$\ict{x,y}f = \frac{-7}6,\;  \ict{y,x}f=\frac{10}3,\;
\ict{x,x+y}f=\frac92,\;\ict{y,x+y}f=0.$$\qed
\end{ex}
One can represent $\cta$ as a sum of iterated constant term
functionals as follows.  Associate to each $\ba\in\aind$ the full flag
of subspaces of $V$:
$$\flag(\ba)=(\cap_{i=1}^n H_{i,\ba}=\{0\},\dots,H_{n-1,\ba}\cap
H_{n,\ba},H_{n,\ba}).$$
Denote $\dim H^n(\compl,\C)$ by $r(\A)$, and define a subset
$\bB\subset\aind$ to be an \emph{orthogonal basis of $\A$} (in degree
$n$) if
\begin{itemize}
\item $|\bB|=r(\A)$
\item for $\ba,\bb\in\bB$ and $\tau\in\perm_n$ the equality
  $\flag(\bb^\tau)=\flag(\ba)$ implies $\bb=\ba$.
\end{itemize}
Then keeping the notation and concepts introduced above, one
obtains the following expression for the constant term.
\begin{prop}
  \label{thm:ires}  
For any orthogonal basis $\bB$ of $\A$ and function $f\in\rat$
one has 
\[ \cta(f) = \sum_{\ba\in\bB}\ict \ba f. \]
\end{prop}
\begin{ex}\label{ex:2dim}
Let $\A=(H_1,\dots,\H_N)$ be a HPA in dimension 2.  Then 
\[ \{(H_1,H_i)|\;i=2,\dots, N\}\subset\A^2\]
is an orthogonal basis.\qed
\end{ex}
Note that the constant term of a function depends on the
hyperplane arrangement. For example, in the notation from above,
$\cta(1)=r(\A)$.
\begin{ex} We compute the constant term of the   function $f=
  \frac{e^{x-2y}}{xy(x+y)}$ with respect to the hyperplane
  arrangement $\A=\{x,y,x+y\}$ using two different orthogonal
  bases:
\[ \cta f = \ict{x,y} f+
\ict{x,x+y}f =\ict{y,x}f+\ict{y,x+y} = \frac{10}3.\] \qed
\end{ex}
To have an efficient method of computation of the constant term,
we need to find such orthogonal bases. Fortunately, they turn out 
to be plentiful. The construction described below is a refined
version of the usual \nbc\ (no broken circuit) bases
(cf. \cite{OT}).  Here we will make an essential use of the
ordering $\prec$ introduced at the beginning of this section.
\begin{prop}[\cite{bern}]\label{thm:nbc}
The set
\begin{multline} \label{eq:def-nbc}
\bcb(\A,\prec) = \{\ba\in\aind|\, \ba\oset \A\text{ and }\\
\text{ for any } H\notin \ba \text{ the set } \{H\}\cup\{G\in \ba|\,
H\prec G\}\text{ is independent}\},
\end{multline} 
is an orthogonal basis of $\A$.\end{prop}
\begin{rem}
1. Note that we are only considering \nbc-bases in degree $n$, which
is the top degree, even though they make sense in all degrees. \\
2. There are examples of orthogonal bases which are not \nbc-bases for
any ordering.  \cite[5.3]{bern}.\\ 
3. The orthogonal basis in Example \ref{ex:2dim} is an \nbc-basis.
\end{rem}

\subsection{Rational sums}
\label{sec:rat}
 Let $V$ be an $n$-dimensional complex vector space and
$\Gamma\subset V$ a lattice of rank $n$. The lattice spans an
$n$-dimensional real subspace $\vr\subset V$ with dual $\vdr\subset V^*$.
Denote the dual lattice of $\Gamma$ by
$\Gamma^*=\Hom(\Gamma,\Z)\subset V^*$ and the group of characters of
$\Gamma$ by $\Gc=\Hom(\Gamma,U(1))\simeq V^*/\Gamma^*$. This last
isomorphism is given by the correspondence $\tau \mapsto t+\Gamma^*$,
with $\tau=\bet|_{\Gamma}$, where  $\bet(v)=e^{\ala
  t(v)}$ for $v\in V$.  

We will say that a central HPA $\A$ in $V$ is \emph{compatible} with
$\Gamma$, or that $(\A,\Gamma)$ is a compatible pair in $V$, if each
hyperplane $H_i$ in the arrangement is the zero-set of some linear
form $x_i\in \Gamma^*$. Define $t\in\vdr$ to be
$\Gamma$-\emph{special} with respect to $\A$ if $t =
\lambda+\sum_{i=1}^N \nu_i x_i$, where $\lambda\in\Gamma^*$,
$\nu_i\in\R$,  and at
most $n-1$ of the coefficients $\{\nu_i\}$ are nonzero. While this
property depends on both $\Gamma$ and $\A$, the reference to
these two objects will be omitted from the notation whenever this causes no
confusion. The set of nonspecial elements is a $\Gamma^*$-invariant
union of open polyhedral chambers.  Since being special is a
$\Gamma^*$-invariant property, we have a well-defined notion of a
\emph{special character} $\tau\in\Gc$, as well.

For an  $n$-tuple $\by=(y_1,\dots,y_n)$ of linear forms
introduce the notation
\[  \Box(\by) = \left\{t\in V^*|\, t=\sum_{i=1}^n \nu_i y_i\text{ with
  }0<\nu_i<1,\, i=1,\dots,n\right\}. \] Assuming that $\by\subset\Gamma^*$,
denote by $\mathrm{vol}_\Gamma(\by)$ the $\Gamma^*$-volume of
$\Box(\by)$, which is always a positive integer. Then we have
\begin{prop} \label{thm:def-tq}
  Let $(\A,\Gamma)$ be a compatible pair in V, and fix $\ba\in\aind$
  and a nonspecial character $\tau\in\Gc$.  Pick a set
  of representative forms $\widehat \ba\subset\Gamma^*$ for $\ba$. \\
  \emph{1.} Then the function
\begin{equation} \label{eq:def-todd}
\td(\Gamma,\ba,\tau) = \frac{\sum \{e_{\tilde t}|\; {e_{\tilde
  t}}_{|\Gamma}=\tau,\,\tilde t\in\Box(\widehat
  \ba)\}}{\mathrm{vol}_\Gamma(\widehat \ba)} \prod_{i=1}^n
  \frac{\ala{x}_{i,\widehat \ba}}{{e}_{i,\widehat \ba}-1}
\end{equation}
in $\mer$ is independent of the choice of forms $\widehat \ba$.\\ 
\emph{2.} The correspondence
\begin{equation}
\label{eq:def-q}
\alpha_\ba \mapsto \td(\Gamma,\ba,\tau) \,\alpha_\ba
\end{equation} induces a  well-defined injection
\[ \iota_\tau^\Gamma:\Omega^n_\A\longrightarrow\Omega_\mathrm{loc}^n(\compl),\]
where $\Omega^n_\A$ is the image of $q$ in degree $n$, and
$\Omega_\mathrm{loc}^n(\compl)$ is the space of holomorphic $n$-forms
defined on the intersection of $\compl$ with a neighborhood of $0$.
\end{prop}
\begin{rem} The  number of terms in the sum in \eqref{eq:def-todd} is
  $\mathrm{vol}_\Gamma(\widehat \ba)$.
\end{rem}
This result allows us to define a deformation 
$\tq:H^n(\compl,\C)\rightarrow\Omega_{\mathrm{loc}}^n(\compl)$ of the
map $q$  via $\tq=\iota_\tau^\Gamma\circ q$. We
emphasize that, while $q$ is defined in all degrees, the deformation
$\tq$ can only be defined in degree $n$.
\begin{cor} For a nonspecial $\tau\in\Gc$, the
functional 
\begin{equation}
  \label{eq:ctatilde}
\widetilde{\mathrm{CT}}_\tau^{\A\Gamma}:\mer\rightarrow\C;\quad
\widetilde{\mathrm{CT}}_\tau^{\A\Gamma}(f)= \sum_i \int_{Z_i}f\tq(b^i)   
\end{equation}
depends only on $\A,\Gamma$ and $\tau$. Its value may be computed via
the formula
\begin{equation}  \label{eq:alg-form}
 \widetilde{\mathrm{CT}}_\tau^{\A\Gamma}(f) = 
\sum_{\ba\in\bB} \ict\ba \left(\td(\Gamma,\ba,\tau)f\right),
\end{equation}
where $\bB$ is an orthogonal basis of $\A$.
\end{cor}
   After these preparations we may formulate the main result:
\begin{thm}[\cite{bern}]
  \label{thm:main-rational} 
  Given a compatible pair $(\A,\Gamma)$, let $f\in\rat$ and
  $\tau={\bet}_{|\Gamma}$. Then the Fourier series
\begin{equation}\label{eq:ber-def}
 B_f^{\A\Gamma}(t) = \sum_{\gamma \in \Gamma\cap\compl} \tau(\gamma)
f(\gamma) \end{equation}
defines a $\Gamma^*$-invariant distribution on $\vdr$ in
the variable $t$, which restricts to a polynomial function on each chamber of
nonspecial elements. Moreover, the formula
\begin{equation}
  \label{eq:main-rational}
  B_f^{\A\Gamma}(t) = (-1)^n\widetilde{\mathrm{CT}}_{\tau}^{\A\Gamma}(f)
\end{equation}
holds for each nonspecial $\tau$.
  \end{thm}
\begin{ex}
An example of the computation of this sum is
\begin{multline}
\label{eq:answer}
B(u,v) = \sum
\frac{e^{u\um+v\un}}{\um^2\un^2(\um+\un)^2},\;m,n,m+n\in\Z^{\neq0} = 
\\ \ict{x,y}\frac{xy}{(1-e^x)(1-e^y)}
\frac{e^{\{ u\} x+\{ v\} y}}{x^2y^2(x+y)^2}
+\ict{x,x+y}\frac{x(x+y)}{(1-e^x)(1-e^{x+y})}
\frac{e^{\{ u-v\} x+\{ v\} (x+y)}}{x^2y^2(x+y)^2} \\
=\ict{x,y}\frac{xy}{(1-e^x)(1-e^y)}
\frac{e^{\{ u\} x+\{ v\} y}+
e^{(1-\{ u-v\}) x+\{ v\} y}}{x^2y^2(x+y)^2}.
\end{multline}
\end{ex}
Here underlining means multiplication by $2\pi \sqrt{-1}$, and $\{u\}$
is the fractional part of $u$. In the course of the computation, we first
rescaled all variables by $2\pi\sqrt{-1}$, then performed the change
of variables 
\[ \{x\rightarrow -x,\,y\rightarrow x+y,\,x+y\rightarrow y\}\]
in the second iterated constant term. Note that the constant term with
respect to $x$ is the same thing as the constant term with respect to
$-x$. The result is a piecewise polynomial function of degree 6 in the
variables $u$ and $v$, with rational coefficients. The complete answer
is to long to write down, but one has, for example,
$B(\frac12,\frac13)=-197/39191040$.  Looking at the answer
\eqref{eq:answer}, it would seem that, say, $B(0,0)$ is not
well-defined, since the fractional part function $u\rightarrow\{u\}$
has a discontinuity at 0. It is clear from the Fourier series,
however, that the function $B(u,v)$ is continuous. This simply means
that the piecewise polynomial functions that we obtain
``miraculously'' agree on the common boundaries of their respective
domains. In particular, setting $u=v=0$ in one of these functions yields
$B(0,0)=-1/30240$, whenever the domain of definition of the function
contains the origin in its closure.\qed

Now we recall the basic steps of the proof of
\thmref{thm:main-rational} given in \cite{bern}. This will be
instructive for the proof of our main result, \thmref{thm:main-trig}.
We refer the reader to the original paper \cite{bern} for details.

\textbf{S1.} Prove the case $n=1$; this follows
from the residue theorem in the complex plane.

\textbf{S2.} Prove the case $n=N$. (Recall that $|\A|=N$.) This can be
treated by considering the product case and then passing to a finite
extension of $\Gamma$ if necessary.

\textbf{S3.} Show that if $N>n$, then the space of rational functions
$\rat$ is spanned by subspaces $R_{\A'}$, where $\A'$ is a nontrivial
subset of $\A$.  This is a variant of a {\em partial fraction
  decomposition} principle in several dimensions.  Much stronger
statements are true, cf.  \cite{zel,bv}. Below we sketch the proof of
a version, which will be used later in the paper.
\begin{prop} \label{thm:rat-parfrac} Let $(\A,\prec)$ be an ordered
essential central HPA, and let $\fA$ be a set of representative linear
forms. Then the space of $\A$-rational functions $\rat$ is linearly
spanned by functions of the form
\begin{equation}  \label{eq:basis}
\frac g{x_{1,\ba}^{\alpha_1}x_{2,\ba}^{\alpha_2}\dots
  x_{n,\ba}^{\alpha_n}},
\end{equation}
where
$\ba\in\bcb(\A,\prec);\,\alpha_1,\alpha_2,\dots\alpha_n\in\Z^{\geq0}$,
and $g$ is a polynomial in the variables $\{x_{j,\ba}|\alpha_j=0\}$.
  \end{prop}
\begin{rem}
  Note that if each of the exponents $\alpha_j$, $j=1,\dots,n$ in
  \eqref{eq:basis} is nonzero, then the only permissible numerator $g$
  is a constant, which can be set to 1. The set of such functions is
  linearly independent and spans a vector space $G_\A$, which is
  independent of the ordering.  This space was introduced by Brion and
  Vergne in \cite{BV}.
\end{rem}
\emph{Proof}: Order all subsequences of the
sequence of forms $(x_1,\dots,x_N)$ lexicographically as follows.
\hspace{1mm} Given two subsequences, \\
$ (y_1,\dots,y_m),(y'_1,\dots,y'_l)\oset\A$,
we will write $(y'_1,\dots,y'_l)\prec (y_1,\dots,y_m)$ if
\begin{itemize}
\item  $l>m$ and $y_{m-k}=y'_{l-k}$ for $k=,\dots,m-1$; or
\item $y'_{l-k} \prec y_{m-k}$, where $k=\min\{k'|\;y_{m-k'}\neq y'_{l-k'}\}$.
\end{itemize}

Now we present an algorithm for exhibiting a rational function
$f\in\rat$ as a linear combination of fractions of the form described
in the Proposition. Recall that a subsequence $(y_1\dots,y_m)\oset \A$
is a \emph{broken circuit} if there exists a minimally linearly
dependent subsequence of the form $(y_0,y_1,\dots,y_m)\oset\A$. In
other words,  a subset of an ordered set of forms is a broken
circuit if and only if there exists a form, preceding all the forms in the
subset, which can be uniquely expressed as a linear combination of the
forms in the subset.  In particular, the forms in the subset itself
should be linearly independent.

Without loss of generality, we can assume that $f$ is a single rational
fraction of the form
\[ f=\frac{g}{\displaystyle \prod_{x\in\fA}x^{\varepsilon(x)}},\]
where $g$ as a polynomial and
$\varepsilon:\fA\rightarrow\Z^{\geq0}$.  Assuming there are
broken circuits in the denominator of $f$, i.e. among the forms
$\{x|\,\varepsilon(x)\neq0\}$, denote by
$\mathrm{bc}(\varepsilon)=(y_1,\dots,y_m)$ the greatest of these
with respect to the lexicographic ordering introduced above.
Also, let $m_\varepsilon=\min\{\varepsilon(y_i), \,i=1,\dots,m\}$
for this broken circuit, and let $y_0\in\fA$ be any element such
that $y_0\prec y_1$ and $(y_0,y_1,\dots,y_m)$ is linearly
dependent.

Then using the obvious identity
\begin{equation}
\label{eq:obv-rat}
\frac{1}{y_1\dots y_m}=
-\sum_{i=1}^m\frac{\lambda_i/\lambda_0}
{\displaystyle \prod_{j\in\overline{0,m}\setminus\{i\}} y_j} 
\quad\mathrm{if}\quad \sum_{i=0}^m\lambda_iy_i=0,
\end{equation}
one can express $f$ as a linear combination of other fractions.
It is easy to check that if $\varepsilon'$ is the exponent
function of one of the new fractions, then
$\mathrm{bc}(\varepsilon')\preceq \mathrm{bc}(\varepsilon)$, and
if $\mathrm{bc}(\varepsilon')= \mathrm{bc}(\varepsilon)$, then
$m_{\varepsilon'}<m_\varepsilon$. We can iterate this elementary
step of the algorithm, applying it to each of the new terms
separately. At the end, we arrive at an expression of $f$ as a
linear combination of fractions, each of which has no broken
circuits among the forms in its denominator. After appropriate
simplifications between the numerators and denominators, one
arrives at the statement of the Proposition. \qed

\textbf{S4.} The proof of \eqref{eq:main-rational} proceeds via an
inductive comparison of the two sides. The induction is carried out on
$|\A|=N$, the number of hyperplanes, and starts with \textbf{S2}.
Thus assume that $N>n$ and that \eqref{eq:main-rational} holds for all
arrangements $\B$ with $|\B|<N$. According to \textbf{S3}, it is
sufficient to consider functions $f$ which are generically regular
along some hyperplane $H\in\A$. Introduce the notation $\A\backslash
H=\A\backslash\{H\}$, and $\arh=\{H\cap L|\,L\in\A\}$, the restriction
of $\A$ onto $H$. Then $f\in R_{\amh}$, and clearly, we have
\begin{equation} \label{eq:bern-ind}
 B_f^{\A\Gamma}(t) = B_f^{\amh\,\Gamma}(t)
-B_{f_{|H}}^{\arh\,\Gamma\cap H}(t_{|H}),
\end{equation}
where $f_{|H}$ and $t_{|H}$ are the natural restrictions. Note that
the compatibility of $\A$ and $\Gamma$ guarantees that $\Gamma\cap H$
is a lattice of full rank in $H$.

Then the proof of the inductive step follows from the relation of the
triple for noncommutative no-broken-circuit bases \cite[Proposition
3.9]{bern}.  This relation is a simple generalization of the
contraction-restriction relation of B\marginpar{who}.  Assume that the
ordering $\prec$ of $\A$, on which so far we have not imposed any
conditions, satisfies the following properties:
\begin{itemize}
\item $H$ is the maximal element of $\A$ with respect to $\prec$,
i.e. $H_N=H$;
\item if $(K,L,M)\oset\A$ and $K\cap H=M\cap H$ then $K\cap H=L\cap
H$.
\end{itemize}
These properties assure that $\prec$ induces a natural ordering
$\prec_{|H}$ on $\arh$.  Now by associating to each $H'\in\arh$ the
{\em minimal} hyperplane $K\in\A$ such that $K\cap H=H'$, we can define a
section $s:\arh\rightarrow\A$ of the canonical map
$\A\rightarrow\arh$. Let $s_H:\arh^{n-1}\rightarrow\A^n$ be the map
induced by $s$ on $\arh^{n-1}$, composed with appending $H$ at the end
of the resulting sequence from $\A^{n-1}$.
Then the relation of the triple reads:
\begin{equation}\label{eq:nbc-ind}
 \bcb(\A,\prec) = \bcb(\amh,\prec)\cup s_H(\bcb(\arh,\prec_{|H})). 
\end{equation}
The comparison of \eqref{eq:bern-ind} and \eqref{eq:nbc-ind} implies
the inductive step. This concludes the proof of the Theorem.\qed

Note the key idea of the proof: we can use the
flexibility provided by Propositions \ref{thm:ires} and \ref{thm:nbc} to
choose a convenient ordering {\em depending} on the function $f$. This
ordering then provides us with a simple expression for the right hand
side of \eqref{eq:main-rational}.

\subsection{Bernoulli Polynomials}\label{sec:polynomial} The function
$B_f^{\A\Gamma}(t)$ introduced in \thmref{thm:main-rational} is
manifestly $\Gamma^*$-invariant, thus it cannot be a polynomial unless
it is constant. The Theorem states that it coincides with a polynomial
when restricted to a chamber of nonspecial elements; this
polynomial will vary from chamber to chamber, however. 
The resulting functions were termed {\em multiple Bernoulli polynomials} in
\cite{bern}. 

Let $\A,\Gamma,f$ be as in the Theorem. For a nonspecial $u\in\vdr$
denote by $\BB^{\A\Gamma}_{[u]f}(t)$ the polynomial function which
coincides with the restriction of $B^{\A\Gamma}_{f}(t)$ to the unique
chamber to which $u$ belongs.  Our immediate goal is to write down a
version of \eqref{eq:main-rational} with the LHS replaced by
$\BB^{\A\Gamma}_{[u]f}(t)$. To see what needs to be changed on the RHS,
recall that the essential ingredient in the definition of
$\ctt_{\tau}^{\A\Gamma}$ was the Todd function \eqref{eq:def-todd}.

Set $t=u$ and rewrite the exponential sum in the numerator as
\[ \sum\{e_{u+w}|\;w\in\Gamma^*,\;u+w\in\Box(\bq)\}.\]
Then it is clear that in order to have a global polynomial such as
$\BB^{\A\Gamma}_{[u]f}(t)$, we  need to replace this sum by
\[ \sum\{e_{t+w}|\;w\in\Gamma^*,\;u+w\in\Box(\bq)\}.\]

It will be convenient to write this modified definition of the Todd
function in terms of the parameter $\mu=t-u$, the {\em shift vector}:
\begin{equation}
  \label{eq:deftodmod}
\td_\mu(\Gamma,\ba,\tau) =
  \frac{\sum \{e_{\tilde t}|\;  
  {e_{\tilde t}}_{|\Gamma}=\tau,\,\tilde t-\mu\in\Box(\widehat
\ba)\}}{\mathrm{vol}_\Gamma(\widehat \ba)}\prod_{i=1}^n 
  \frac{\ala{x}_{i,\bah}}{\mathbf{e}_{i,\bah}-1}.
\end{equation}
It is easy to check that this definition gives rise to a
consistent deformation $\tilde q_{\mu,\tau}^\Gamma$ of $q$ and to
a deformed constant term functional $\ctt_{\mu,\tau}^{\A\Gamma}$.
Then the new variant of \thmref{thm:main-rational} reads:
\begin{thm}\label{thm:var-main-rational}
  Let $\A,\Gamma$ be compatible, and $f\in\rat$. Let $u,t\in\vdr$,
  with $u$ nonspecial, and set $\mu=t-u$. Then
  \begin{equation}
    \label{eq:var-main-rational}
    \BB^{\A\Gamma}_{[u]f}(t)=(-1)^n\ctt_{\mu,\tau}^{\A\Gamma}(f).
  \end{equation}
\end{thm}
An advantage of this formulation is that it allows us to evaluate the function
$B^{\A\Gamma}_{f}$ at special values of $t$ as well.
\begin{cor} \label{thm:rat-gen}
Assume that $\A$ and $\Gamma$ are compatible and
$f\in\rat$ is such that the series \eqref{eq:ber-def} defining
$B^{\A\Gamma}_{f}$ is absolutely convergent. Then for an arbitrary,
possibly special $t\in\vdr$ one has 
\begin{equation}
  \label{eq:ext-rational-sing}
  B_f^{\A\Gamma}(t) =
  (-1)^n\widetilde{\mathrm{CT}}_{\mu,\tau}^{\A\Gamma}(f),
\end{equation}
where, as usual, $\tau=e_t|_{\Gamma}$, and $\mu$ is a sufficiently
small vector in $\vdr$ such that $t-\mu$ is not special.
\end{cor}
When $\mu=0$, then this statement coincides with the statement in
\thmref{thm:main-rational}. We will give an example of a calculation
where $t$ is special in the next section.

\section{Affine arrangements}\label{sec:affine}

Here we describe a generalization of the results of the
previous section to affine arrangements. 

An affine hyperplane is  one that does not necessarily go through the
origin. For each such hyperplane $H$, there is a unique parallel
hyperplane $H^\circ$ containing the origin. Similarly, for each affine
linear form $x$ there is a unique linear form $x^\circ$ such that
$x-x^\circ$ is a constant.

Let $\A$ be a collection of affine hyperplanes in the $n$-dimensional
vector space $V$. Consistently with the notation introduced above, we
have a central HPA $\Ac$ consisting of translates of the elements
of $\A$, and for each $n$-tuple $\ba\in\A^n$ we have a corresponding
$\ba^\circ\in\A^{\circ n}$. We will use a similar notation for
$n$-tuples of forms. Note that it is possible that $|\A|>|\Ac|$, since
we allow parallel planes.

We will call an affine HPA $\A$ and a lattice $\Gamma$ {\em
  compatible} if $\Ac$ and $\Gamma$ are. In this situation we will
assume that the chosen representative forms satisfy
$\fA^\circ\subset\Gs$, and for $y\in\fA$ one has
$y^\circ\in\fA^\circ$.

 Define $p\in V$ to be a \emph{vertex} of $\A$ if $\cap\A_p=\{p\}$,
where $\A_p=\{H\in\A|\,p\in H\}$. We assume that the set $\vert(\A)$ of
vertices of $\A$ is nonempty.  Note that it is not assumed that the
constants $x-x^\circ$ are real, thus it is possible that $p\notin \vr$.

The generalization now goes as follows. As each arrangement $\A_p$ is
essentially a central arrangement, it has its own constant term
functional $\ct^{\A_p}$ taken at the point $p\in V$. Then we define
the constant term of $\A$ by
\begin{equation}
  \label{eq:ct=ctp}
  \cta = \sum_{p\in \vert(\A)} {\ct}^{\A_p},
\end{equation}
The definition of the deformation of $\ct^{\A_p}$ is similar to the
deformation in the central case. One needs to generalize the
definition of \eqref{eq:deftodmod} for an $n$-tuple $\ba$ of {\em
  affine} linear forms as follows:
\begin{equation} \label{eq:def-aff-todd}
\td_\mu(\Gamma,\ba,\tau) =
  \frac{\sum \{e_{\tilde t}|\;  
  {e_{\tilde t}}_{|\Gamma}=\tau,\,\tilde t-\mu\in\Box(\widehat
  \ba^\circ)\}}{\mathrm{vol}_\Gamma(\widehat \ba^\circ)} 
  \prod_{i=1}^n \frac{\ala{x}_{i,\widehat \ba}}{{e}_{i,\widehat
      \ba^\circ}-1},  
\end{equation}
The statements of \propref{thm:def-tq} still hold, since the affine
forms in $\fA_p$ satisfy exactly the same linear relations as the ones
in $\fA^\circ_p$.  
Merging \eqref{eq:ct=ctp} with the results of the previous section, 
\eqref{eq:alg-form}, \eqref{eq:deftodmod},
\eqref{eq:ext-rational-sing}, etc., we can write  down the analogous
equalities in the affine case:
\begin{eqnarray}  \label{eq:aff-alg-form}\displaystyle
 \widetilde{\mathrm{CT}}_{\mu,\tau}^{\A\Gamma}(f) = 
\sum_{p\in \vert(\A)}\widetilde{\mathrm{CT}}_{\mu,\tau}^{\A_p\Gamma}(f),\\
\label{eq:def-mu-p}
\widetilde{\mathrm{CT}}_{\mu,\tau}^{\A_p\Gamma}(f)=
\sum_{\ba\in\bB_p} \ict\ba
\left(\td_\mu(\Gamma,\ba,\tau)f\right),
\end{eqnarray}
where $\bB_p$ is an orthogonal basis of $\A_p$ consisting of affine
linear forms. 

Under the same conditions as in \corref{thm:rat-gen}, and using
\eqref{eq:aff-alg-form} and \eqref{eq:def-mu-p}, we have
\begin{equation}
  \label{eq:main-affine}
B_f^{\A\Gamma}(t) = (-1)^n  
\widetilde{\mathrm{CT}}_{\mu,\tau}^{\A\Gamma}(f),
\end{equation}
As we will prove a more complicated version of this statement in the
trigonometric case, the proof of this equality will be omitted. 

\begin{rem}
  1.  Note that this statement is a common generalization of the
  results of \cite{bern}, sketched in the previous section, and
  the formulae of Brion-Vergne given in \cite{BV}, where, in
  particular, the case of
  a single vertex at a generic point was considered.\\
  2. It is very important that in \eqref{eq:deftodmod} the {\em
    affine} forms $\bq$ appear in  $x_{i,\bq}$ instead of the corresponding
  linear forms, as is the case in the rest of the formula.  As an
  exercise, the reader may check that if $\A=\A_p$ has a single
  vertex $p\notin\Gamma$, then
  $\widetilde{\ct}_{\mu,\tau}^{\A_p\Gamma}(1)=0$.
\end{rem}

We end this short section with an example.
\begin{ex}
We compute the sum
\begin{equation}
  \label{eq:aff-example}
\sum \frac1{mn(2m+n-1)},\quad
m,n\in\Z^{\neq0}, \, 2m+n\neq1.
\end{equation}
Here $\fA=\{x,y,2x+y-1\}$, $\Gamma=\Z^2$ and $\tau$ is the trivial
character. This series converges absolutely, albeit painfully slowly.
The arrangement has 3 simple vertices, each one contributing a single
constant term. Since $\tau$ is singular, we need to choose a small
shift vector $\mu$. We choose $\mu=-\epsilon y-\epsilon^2x$, where
$\epsilon$ is a small positive number. If we choose $\mu$ differently,
the formulas change somewhat, but the result, naturally, remains the
same.

Again, rescaling the variables by $2\pi\sqrt{-1}$ and dividing the sum
by $(2\pi\sqrt{-1})^3$, we arrive at the
result
\begin{multline} \ict{x,\,y}\frac1{\tt x\tt y(2x+y-\ala1)}\\
+\ict{x,\,2x+y-\ala1}\frac1{\tt x\tt{2x+y} y}\\
+\ict{2x+y-\ala1,\,y}  \frac12\frac{\left(e^y+e^{x+y}\right)}
  {\tt y\tt{2x+y} x}.
\end{multline}
As all threes vertices are simple, i.e. each is contained in exactly
two hyperplanes, one can replace the iterated constant terms with the
ordinary constant term. After the appropriate substitutions we obtain
\[
\mathrm{CT}\; \frac1{\tt x\tt y}
\left(\frac1{2x+y-\ala1}+\frac1{y-2x+\ala1}
+\frac{e^y-e^{\frac{x+y}2}}{x-y+\ala1}\right)=\frac{\pi^2-8}{(\ul)^3}.
\]
Thus the value of our infinite sum \eqref{eq:aff-example} is
$\pi^2-8$.\qed
\end{ex}
While it is possible to compute this answer by some ad hoc method as
well, the advantage of our formula is that it does not become more
complicated as the powers of the linear forms in the denominator of
the function $f$ increase.

\section{The trigonometric case}

\subsection{Toric arrangements and rational trigonometric sums} In this section we present a periodic version of the theory, i.e. when
hyperplanes are replaced by hypertori on a torus.  We keep the
notation of the previous section: $\A$ is an essential affine HPA and
$\Ac$ is the associated central arrangement in a complex vector space
$V$. Let $\Theta\subset V$ be a compatible rank-$n$ lattice and denote
by $T$ the complexified torus $V/\Theta$. Recall that $\Theta$ defines
a real subspace $\vr\subset V$.

In view of the compatibility condition, each $H\in\A$ defines a
hypertorus $H+\Theta\subset T$; this results in an arrangement
$\A/\Theta$ of hypertori in $T$.  It will be convenient to work with
the inverse image of this arrangement under the natural map
$V\rightarrow T$. This inverse image is an infinite, periodic HPA on
V, which we denote by $\At$.  Recall that $\Atp$ stands for the
arrangement of those hyperplanes from $\At$ which go through the point
$p\in V$.  The arrangement $\At$ has infinitely many vertices, but we
are interested in them up to a translation by a vector from $\Theta$
only; define $\vert(\A/\Theta)=\vert(\At)/\Theta$. Then elements of
$\vert(\A/\Theta)$ are the vertices of the toric arrangement
$\A/\Theta$.

Fix an ordering $\prec$ on $\Ac$, choose a set $\Ahc$ of
representative forms which are {\em minimal} elements of $\Theta^*$,
and let $\fA$ be the corresponding the of affine linear forms.

Now we introduce the periodic version of $\A$-rational functions.  Let
$\ctc$ be the polynomial ring generated by the functions
$\be_y$, $y\in\Theta^*$, and let $\catc$ be the same ring,
with the functions $1-\bey$ inverted whenever $y\in\widehat\A$. We
will think of these spaces as of spaces of functions on $T$ or
periodic functions on $V$, interchangeably. The elements of $\catc$
will be called \emph{trigonometric $(\A,\Theta)$-rational functions}.

Note that the constant term functional $\ct^{\Atp}$ restricted to
$\catc$ is independent of the choice of the vertex $p\in\vert(\At)$
modulo $\Theta$.  Then, just as in the previous section, we may define a
constant term functional $\ct^{\A/\Theta}:\catc\rightarrow\C $ given
by the finite sum
\[ \ct^{\A/\Theta} = \sum_{p\in \vert(\A/\Theta)} {\ct}^{\Atp}.\]
Here we took up the somewhat sloppy convention that summation over
$p\in \vert(\A/\Theta)$ means taking a representative $p$ from each of
the $\Theta$-equivalent classes of vertices of $\At$.

Now let $\Gamma$ be a lattice containing $\Theta$, and let 
$f\in\catc$. Fix a character $\tau\in\Hom(\Gamma/\Theta,U(1))$; this
may be thought of as an element of $\Gc$ which is trivial on $\Theta$.
We are interested in rational trigonometric sums of the form
\begin{equation} \label{vdef}
Z_f^{\A\Gamma/\Theta}(\tau) =\sum_{\gamma\in(\Gamma\cap U(\At))/\Theta}
\tau(\gamma)f(\gamma).
\end{equation}
The interest in such sums was sparked by a formula given by
E. Verlinde for the dimension of conformal blocks of the WZW theory
\cite{Ver}. We will study this formula later in the paper.
We would like to write down a localized formula for the sums
\eqref{vdef} similar to \eqref{eq:main-rational}.  Recall that,
originally, when we treated the rational case in \secref{sec:rat}, we
excluded the special characters. This was partly justified because the
Fourier series $B^{\A\Gamma}_f(t)$ may very well have singularities
for special values of $t$. Later we treated the case of special values
in \secref{sec:polynomial}.

In the trigonometric case, which we are investigating in this section,
we cannot exclude the special characters  since every sum is
meaningful. Here, however, there is no applicable notion of
continuity, such as the one used in \secref{sec:polynomial}.
 
 It turns out that the shift vector $\mu$, which played an auxiliary
 role in \secref{sec:polynomial}, in the trigonometric case becomes
 essential. Going back to \eqref{eq:def-aff-todd} and
 \eqref{eq:aff-alg-form}, we can write down a natural definition for
 the deformed constant term of the toric arrangement $\A/\Theta$ with a
 $\mu$-shift:
\[ \widetilde{\ct}^{\A\Gamma/\Theta}_{\mu,\tau} =
\sum_{p\in
  \vert(\A/\Theta)}\widetilde{\ct}^{\Atp\Gamma}_{\mu,\tau}, \] where
the functional $\widetilde{\ct}^{\Atp\Gamma}_{\mu,\tau}$ is defined in
\eqref{eq:def-mu-p}, but here we consider it to be restricted onto
$\catc$.

Turning to the function $f\in\catc$, define
\begin{multline}
  \label{eq:def-delta}
 \Delta_f = \{ \mu\in\vdr|\; \textrm{for }u\in\vr\textrm{ not a pole
 of }f\text{ and for all nonzero }v\in\vr\\
\lim_{s\in\R,\,s\rightarrow +\infty}\be_\mu(\sqrt{-1}(u+sv))f(\sqrt{-1}(u+sv)) = 0\}. 
\end{multline}
Note that this definition is independent of the HPA, and depends on
the function only.
We introduce the associated linear spaces
\[ \catcm = \{f\in\catc|\;\mu\in\Delta_f\}.\]

By definition, after putting its terms over a common denominator,
any $f\in\catc$ may be represented in the form
\begin{equation}
  \label{eq:fform}
  f=\frac{\sum_{w\in\Theta^*}\lambda_w\bew}
{\prod_{y\in\fA}(1-{e}_{y})^{\varepsilon(y)}} ,
\end{equation}
where $\varepsilon:\fA\rightarrow \Z^{\geq0}$ and the complex numbers
$\lambda_w$ vanish for all but finitely many $w\in\Theta^*$.
Introduce the finite set
$N_{f,\varepsilon}=\{w\in\Theta^*|\,\lambda_w\neq0\}$ and the convex
polytope $D_{\varepsilon}=\{\sum_{y\in\fA} \nu(y)\varepsilon(y)
y^\circ|\,0<\nu(y)<1\}$. We will use the convention that
$\bD_{\varepsilon}$ is the closure of $D_{\varepsilon}$ unless
$\varepsilon$ is $0$, i.e. when there is nothing in the denominator.
Then $D_{\varepsilon}=\emptyset$ and $\bD_{\varepsilon}=\{0\}$.
\begin{prop}
  \label{thm:dv}
Given a function $f\in\catc$ in the
  form \eqref{eq:fform}, one has $\mu\in\Delta_f$ if and only if
  $D_{\varepsilon}$ is a nonempty open subset of $\vdr$, and $\mu+w\in
  D_{\varepsilon}$ for every $w\in N_{f,\varepsilon}$.
\end{prop}
 Note that the set $D_{\varepsilon}$ is {\em nonempty and open} exactly
when the linear forms which actually appear in the denominator span
$\vdr$. To put the statement in formulas, introduce
\begin{equation}
  \label{eq:def-deltazero}
  \Delta_f^0 = \bigcap_{w\in N_{f,\varepsilon}}(D_{\varepsilon}-w)\quad\text{and}\quad
\bar \Delta_f= \bigcap_{w\in N_{f,\varepsilon}}(\bD_{\varepsilon}-w).
\end{equation}
Then the statement is that $\Delta_f=\inter(\Delta_f^0)=\inter (\bar
\Delta_f)$, where $\inter(S)$ stands for the interior of the set $S$. In fact, if
$D_{\varepsilon}$ is nonempty and open, then $\Delta_f=\Delta_f^0$,
and the closure of $\Delta_f$ is $\bar\Delta_f$.

{\em Proof of Lemma}: The ``if'' part is easy, because if the
conditions on $\mu$ in the Proposition hold, then one can represent
each term in $e_\mu f$ as a product of functions of the form
$\be_{\nu(y) y^\circ}/(1-\be_{y})$, where $0<\nu_{y}<1$ and $y\in\fA$.  As
long as the linear forms $\{y^\circ|\,\varepsilon(y)\neq0\}$ do not lie in a
hyperplane in $\vdr$, for any nonzero $v\in\vr$ there will be at least
one form $y$ for which $y^\circ(v)\neq0$. In this situation the
condition in \eqref{eq:def-delta} clearly holds. 

For the ``only if'' part note that if there is a nonzero $v\in\vr$
such that $y^\circ(v)=0$ for each $y$ in the denominator, then for
such $v$ the condition in \eqref{eq:def-delta} cannot hold, thus it is
necessary that $D_{\varepsilon}$ be nonempty and open. We assume
this from now on. 

Now, ad absurdum, suppose that there is a $w\in
N_{f,\varepsilon}$, for which $\mu+w\notin D_{\varepsilon}$.
Since $D_{\varepsilon}$ is convex, there is a hyperplane
separating $\mu+w$ from $D_{\varepsilon}$. In other words, there
is a $v\in\vr$ such that 
\[
\mu(v)+w(v)\geq
\sum_{y\in\fA}\varepsilon(y)(y^\circ(v)+|y^\circ(v)|)/2.
\]
Clearly, when we restrict $f$ to the line $\sqrt{-1}(u+sv)$ for
some $u$ and let $s\rightarrow-\infty$, then the dominant
contribution in \eqref{eq:fform} comes from those terms of the form
$\lambda_w \bew/{\prod_{i=1}^M(1-\mathbf{e}_{y_i})}$ for which
$w(v)$ has the maximal value $m=\max_{w\in
  N_{f,\varepsilon}}w(v)$. Denoting the set of such $w$s by
$M_v(f)=\{w\in N_{f,\varepsilon}|\,w(v)=m\}$, we can write the
dominant contribution as
\[ \frac{e^{-2\pi (m+\mu(v))s}}
{\prod_{y\in\fA}\left(1-a_ye^{-2\pi y^\circ(v)s)}\right)^{\varepsilon(y)}}\sum_{w\in
M_v(f)} \lambda_w \bew(u),
\]
where $a_y$ is a nonzero constant: $a_y=e^{-2\pi y(u)}$.  Note that
for generic $u$, the coefficient in front of the fraction is not 0
because $M_v(f)$ is nonempty. Thus this dominant contribution does not
vanish as $s\rightarrow-\infty$ since $m+\mu(v)\geq
\sum_{y\in\fA}\varepsilon(y)(y^\circ(v)+|y^\circ(v)|)/2$.  This means
that $f$ does not satisfy \eqref{eq:def-delta}, and this contradicts
our assumption. The proof is complete.\qed

\subsection{The main result}
Now we can formulate the main result of the paper:
\begin{thm}
  \label{thm:main-trig}
  Let $\Theta\subset\Gamma$ be lattices compatible with an affine
  HPA $\A$ in an $n$-dimensional complex vector space $V$. Let
  $f\in\catcm$ for some $\mu\in\vdr$. Then for  $t\in\Ts$
  such that $t-\mu$ is not $\Gamma$-special, the identity
\begin{equation}
  \label{eq:main-trig}
Z_f^{\A\Gamma/\Theta}(\tau) =
\widetilde{\ct}^{\A\Gamma/\Theta}_{\mu,\tau}(f)  
\end{equation}
holds, where $\tau=\bet|_{\Gamma/\Theta}$.
\end{thm}
For computational purposes the iterated constant term form of
$\widetilde{\ct}^{\A\Gamma/\Theta}_{\mu,\tau}(f) $  is useful:
\begin{equation}
  \label{eq:itrestrig}
  Z_f^{\A\Gamma/\Theta}(\tau) = \sum_{p\in \vert(\A/\Theta)}
\sum_{\ba\in\bB_p} \ict\ba
\left(\td_\mu(\Gamma,\ba,\tau)f\right), 
\end{equation}
where $\bB_p$ is an orthogonal basis of the arrangement $\Atp$ (see the
definitions at the start of this section). 
\begin{ex}\label{ex:b2}
  The simplest nontrivial example with multiple vertices is provided
  by the arrangement corresponding to the Lie algebra $B_2$.  Set
  $\fA=\{x,y,x+y,x-y\}$, $\Theta=\Z^2$, $\Gamma=\Theta/k$. Introduce
  the shorthand
\[ T(x)={\tt x\tt{-x}},\quad  \delta(x,y)=\frac1{T(x)T(y)T(x+y)T(x-y)}.\]
Let $f(x,y)=\delta(\ala x,\ala y)$ and $\tau=1$. Our aim is to
compute the sum
\[ Z_f^{\A\Gamma/\Theta}(\tau)=\sum f( /ik, /jk),
\quad 0<i\neq j<k,\;i+j\neq k.\]
The vertices of $\A/\Theta$ are at $(0,0)$ and at
$(\frac12,\frac12)$, with corresponding orthogonal bases
$\{(x,y),(x,x+y),(x,x-y)\}$ and $(x+y-1,x-y)$. Now we need to choose a
$\mu\in\Delta_f$. We will take one of the simplest choices:
$\mu=-\epsilon x-\epsilon^2 y$, but note that we might have picked,
for example, $3x-y-\epsilon x-\epsilon^2y$, which would have
resulted in a rather different formula.

Again, we rescale by $2\pi\sqrt{-1}$ and obtain
\begin{multline}
 \ict{x,y}\frac{k^2xy\delta(x,y)}{\tt{kx}\tt{ky}}+
\ict{x,\,x+y}\frac{k^2x(x+y)\delta(x,y)}{\tt{kx}\tt{kx+ky}}\\
+\ict{x,\,x-y}\frac{k^2x(y-x)\delta(x,y)}{\tt{kx}\tt{ky-kx}}
+\ict{x+y-\ala1,\,x-y}\frac{k^2(x+y-\ala1)(x-y)(1+e^{kx})\delta(x,y)}{2\tt{kx+ky}\tt{kx-ky}}.
\end{multline}
After the usual change of variables this equals
\begin{multline}
k^2\ict{x,y}\frac{xy}{\tt{kx}\tt{ky}}
\left(\vphantom{\frac12}
\delta(x,y)+\delta(x,y-x)+\delta(x,x+y)+
\right. \\
\frac{1+(-1)^ke^{\frac{kx+ky}2}}2
\left.\delta\left(\frac{x+y+\ala1}2,\frac{x-y+\ala1}2\right)
\right).
\end{multline}
The answer is the somewhat intimidating
\[\frac{k^8+60k^6+5523k^4+133377/2+(-1)^k(525k^4
+5250k^2+30975/2)}{240\cdot8!},\]
which, nevertheless, does reduce to $0,1/8,8/25,10/9$ for $k=3,4,5,6$,
respectively.\qed  
\end{ex}
 If $\Delta_f=\emptyset$, then the statement of the Theorem is
  vacuous.  One can still use \eqref{eq:main-trig} to compute
  any rational trigonometric sum using the identity
\[ 1= \frac1{1-z}+\frac1{1-z^{-1}}.\]
Indeed, the identity clearly implies
\begin{lem} The space $\catc$ is linearly spanned by the linear spaces
$$\{\catcm|\,\mu\in\vdr\}.$$
\end{lem} 
Thus any $f\in\catc$ with $\Delta_f=0$ can be represented as a sum of
terms, each of which has a nonempty $\Delta$, and then one can apply
the theorem to each term separately. The simplest example of this is
\begin{multline}
\sum_{\omega^k=1,\,\omega\neq 1} 1 = 
\ct\left[\frac{kx}{1-e^{kx}}\frac1{1-e^{-x}}\right]+
\ct\left[\frac{kxe^{kx}}{1-e^{kx}}\frac1{1-e^{x}}\right]=\\
\frac{k-1}2+\frac{k-1}2 = k-1.\end{multline}

\subsection{The proof} The proof of the theorem is parallel to that of
\thmref{thm:main-rational}. \\
\textbf{Step 1.} For the case of $n=1$ we may identify $V$ with $\C$
and $\Theta$ with $\Z$. Denote the coordinate on
$T\simeq\C\backslash\{0\}$ by $z$. The following statement is an
immediate consequence of the residue theorem in $\C$:
\begin{lem}
  \label{thm:one-dim}
Let $F(z)$ be a rational function,
$k\in\Z^{>0}$ and $l\in\Z$. Assume that $l'$ is an integer, such
that $l'\cong l\mod k$ and $z^{l'}F(z)/(1-z^k)$ vanishes both at
$0$ and at $\infty$. Then we have
\[\sum_{\overset{\omega^k=1}{\omega\notin \mathrm{Pole}(F)}} \omega^lF(\omega)=
\sum_{p\in \mathrm{Pole}(F)} \res_{z=p} \frac{dz}z\frac{kz^{ l'}}{1-z^k}F(z),
\]
where $\mathrm{Pole}(F)$ is the set of poles of the function $F$.
\end{lem}
Now the rank-1 case of \eqref{eq:main-rational} easily follows after
performing the change of variables $z\rightarrow e^{\ala x}$.\\
\textbf{Step 2.} Here we consider the case $|\Ac|=n=\dim V$.
Let $\cx_1,\dots,\cx_n\in\Theta^*$ be minimal defining linear forms
for $\Ac$ and let $\beta_1\cx_1,\dots,\beta_n\cx_n$ be the
corresponding minimal elements of $\Gamma^*$. Note that all the
$\beta_i$s are integers since $\Gs\subset\Ts$. Define the lattices
$\Theta_\A$ and $\Gamma_\A$ through their duals:
\[
\Theta^*_\A=\left\{\left.\sum_{i=1}^nk_i\cx_i\right|\;k_i\in\Z\right\},
\quad
\Gamma^*_\A=\left\{\left.\sum_{i=1}^nk_i\beta_i\cx_i\right|\;k_i\in\Z\right\}.\]
Then we have $\Theta\subset\Theta_\A$ and
$\Gamma\subset\Gamma_\A$, but not necessarily
$\Theta_\A\subset\Gamma$.

We prove \eqref{eq:main-trig} in three steps. We fix the data of
the function $f\in\catcm$ and $\tau\in\Hom(\Gamma/\Theta,U(1))$.

One can easily see that the equality
\eqref{eq:main-trig} holds for the pair of lattices
$\Theta_\A\subset\Gamma_\A$, since in this case both sides are simply
products of the 1-dimensional case proved in \textbf{Step
1}. To apply the theorem in this case, we need to assume that $\tau$
is trivial on $\Theta_\A$.

Replacing $\Theta_\A$ with $\Theta$ is painless; if $\tau$ is trivial
on $\Theta_\A$, then both sides of the equality are simply multiplied
by $|\Theta_\A/\Theta|$ because of the $\Theta_\A$-periodicity of the
data. Similarly, it is easy to see that if $\tau$ is nontrivial on
$\Theta_\A$, then both sides of \eqref{eq:main-trig} vanish because
they are multiplied by the sum of the values of a nontrivial
character of $\Theta_\A/\Theta$.

To pass from $\Gamma_\A$ to $\Gamma$, first recall the fact that for a
finite group G
\begin{equation}\label{eq:finite}
\sum_{\chi\in R(G)}\chi(g) = \begin{cases}
|G|,\; \text{if } g=e,\\ 0,\;\text{otherwise},
\end{cases}
\end{equation}
where $R(G)$ is the set of irreducible characters and $e$ is the unit
element. Applying this to the group $\Gamma_\A/\Gamma$, one can see that
\[ Z_f^{\A\Gamma/\Theta}(\tau) =\frac1{|\Gamma_\A/\Gamma|}
\sum_{\tau'\in\Gc/\Gc_\A}Z^{\A\Gamma_\A/\Theta}_f(\tau\tau').  \] 
Checking the $\widetilde{\ct}$ side, we see that the change in
the exponential sum in the definition \eqref{eq:deftodmod}, when
$\Gamma_\A$ is replaced by $\Gamma$, reproduces the same
relation:
\[ \widetilde{\ct}_{\tau,\mu}^{\A\Gamma/\Theta}(f) =\frac1{|\Gamma_\A/\Gamma|}
\sum_{\tau'\in\Gc/\Gc_\A}
\widetilde{\ct}^{\A\Gamma_\A/\Theta}_{\tau\tau',\mu}(f).
\] Thus knowing \eqref{eq:main-trig} for $\Gamma_\A$ implies the same
equality for $\Gamma$.\qed \\ \textbf{Step 3.} While the other steps
are substantially analogous to the rational case, here there are
significant differences.  One problem is that the direct analog of the
partial fractions approximation principle formulated in
\propref{thm:rat-parfrac} does not hold, i.e. the space of functions
which have their poles on essential subsets of $\A/\Theta$ does not
span $\catc$. Also, even if such decomposition $f=\sum g_i$ existed,
we would need to make sure that we have $\mu\in \Delta_{g_i}$ for each
summand $g_i$.

We start with the trigonometric analog of \eqref{eq:obv-rat}.
\begin{lem}
\label{thm:obv-trig}
Let  $y_0,\dots,y_m$ be affine linear forms on $V$, satisfying
$y_0+y_1+\dots+y_m=0$, and such that the
linear forms $\{y_i^\circ\}_{i=0}^m$ are in $\vdr$.  Let 
$\tilde\mu=\sum_{i=1}^m\alpha_iy^\circ_i$ with 
\[ 0<\alpha_1<\dots<\alpha_m<1.\]
Then the identity
\begin{equation}
  \label{eq:obv-trig}
  \sum_{i=0}^m\frac{\prod_{j=0}^{i-1}\be_{y_j}}{\prod_{j\neq i}(1-\be_{y_j})}=0
\end{equation}
holds, and every term in the sum is a trigonometric rational
function which has $\tilde \mu$ in its $\Delta^0$
(cf. \eqref{eq:def-deltazero}). 
\end{lem}
\begin{rem}
We had to allow for the possibility of $m<n$ in this Lemma, hence
the somewhat awkward formulation using $\Delta^0$ instead of
$\Delta$.  
\end{rem}
{\em Proof}. The identity maybe easily checked by multiplying through with
$\prod_{i=0}^m(1-\be_{y_i})$. To check that $\tilde\mu$ is in the $\Delta^0$
of the $i$th term, we express $\tilde\mu+\sum_{j=0}^{i-1}y^\circ_j$ as a
linear combination of the $y^\circ$ variables less $y^\circ_i$:
\[ \tilde\mu +\sum_{j=0}^{i-1}y^\circ_j=(1-\alpha_i)y^\circ_0+
\sum_{j=1}^{i-1}(1-\alpha_i+\alpha_j)y^\circ_j
+\sum_{j=i+1}^{m}(\alpha_j-\alpha_i)y^\circ_j.
\]
Since all the coefficients of this linear combination are between $0$
and $1$, the statement now follows from the definition of $\Delta^0$. \qed

Using the Lemma we can adapt the algorithm in the proof of
\propref{thm:rat-parfrac} to the trigonometric case, as follows.
Take a $\mu$ which is not special with respect to $\Ac$ and
$\Theta$, and fix a function $f\in\catcm$. Observe that  for a
nonspecial $\mu$, the space $\catcm$ does not change as we vary
$\mu$ in a small neighborhood. In other words, we may switch
$\mu$ for a nearby vector if it is necessary.

According to \propref{thm:dv}, $f$ may be represented as a sum of
elements of $\catcm$ of the form
\begin{equation}
  \label{eq:mform}
 \frac{\bet}{\prod_{y\in\fA}(1-\be_{y})^{\varepsilon(y)}},
\end{equation}
where $\varepsilon:\fA\rightarrow \Z^{\geq0}$.  We may assume without
loss of generality that $f$ has this form to begin with.

Following the blueprint of the partial fraction decomposition in
the rational case, we again find an ordered $m$-tuple
$(y_1,\dots,y_m)$ among $\{y|\,\varepsilon(y)\neq0\}$, such that
$(y^\circ_1,\dots,y^\circ_m)$ is the largest possible broken
circuit with respect to the lexicographic ordering. Then, by
definition, there is $y^\circ_0\in\fA^\circ$ such that
$y^\circ_0\prec y^\circ_1$, and a relation
$\sum_{j=0}^m\lambda_jy_j=0$ holds, where $\lambda_j\in\Z$, $j\in
\overline{0,m}$. Now let
\[ \ty_i =\lambda_iy_i,\;i\in\fs m\quad\text{and}\quad
 \ty_0=-\ty_1-\dots-\ty_m.\]
By separating the factors in
 the denominator corresponding to this broken circuit and applying
 the formula for finite geometric progressions, we can write
\begin{equation}
  \label{eq:geomprog}
f=  \frac{\bet}{\prod_{y\in\fA} (1-\be_{y})^{\varepsilon'(y)}} \prod_{j=1}^m
\frac{1+\be_{y_j}+\dots+\be_{\lambda_{j-1}y_j}}{1-\be_{\ty_j}}.
\end{equation}
Expanding the numerator, we obtain a representation of $f$ as a sum of
terms of the form \eqref{eq:mform}, with each term containing $\mu$ in its
$\Delta$. We will consider each of these terms separately.

Using \propref{thm:dv} yet again, we see that $\mu$ may be split
into a sum of two contributions: $\mu=\mu'+\mu''$, with $\mu'$ in
the $\Delta^0$ of the first part of \eqref{eq:geomprog} and
$\mu''$ is in the $\Delta^0$ of the product. Now we use our
freedom of varying $\mu$ in a small neighborhood to make sure
that $\mu''$ is not special with respect to
$(\ty_0^\circ,\dots,\ty_m^\circ)$; this implies that the
coefficients $\alpha_i$ in the representation
$\mu''=\sum_{i=1}^m\alpha_i\ty_i^\circ$ are all different.
Indeed, if, say, we had $\alpha=\alpha_1=\alpha_2$, then we could
represent $\mu''$ as a linear combination of the variables
$\ty_0, \ty_3,\ty_4,\dots,\ty_m$ only, by subtracting
$0=\alpha\sum_{j=0}^m \ty_j$ from the original representation.

Now we renumber the variables $\ty_j$ according to the increasing
order of the coefficients $\alpha_j$, and then apply
\lemref{thm:obv-trig} to $(\ty_0,\ty_1,\dots,\ty_m)$ and
$\tilde\mu=\mu''$. This allows us to use \eqref{eq:obv-trig} to
eliminate the broken circuit.  Again, similarly to the proof of
\propref{thm:rat-parfrac}, we conclude that after performing this
transformation of the representation of $f$, the parameter
\[ \min_{1\leq i\leq m}\left\{\sum\varepsilon(y), \,y^\circ=y_i^\circ\right\}
\]
in the representation \eqref{eq:mform} decreases, and the
lexicographicly largest broken circuit among the forms in the
denominator does not increase.

Note, however, that our manipulations have a price. As the linear
forms $\ty_j$, $j=0,\dots,m$, are not necessarily in $\fA$, the new
fractions are not going to be in $\catcm$. Rather, they will be in
$\C^\mu_\B[\Tc]$ for some affine HPA $\B$ with $\B^\circ=\A^\circ$. 
It is clear that the new affine forms that we need to allow are
translates of $\cy_0,\dots,\cy_m$: for $\cy_j$, where $j\neq0$, these
translates have the form $y_j+i/\lambda_j$, $i=0,\dots,\lambda_j-1$,
but for $\cy_0$ they have an additional shift.

Iterating the procedure, we arrive at the following analog of
\propref{thm:rat-parfrac}:
\begin{prop}
\label{thm:trig-parfrac}
Let $\A$ be an affine HPA in $V$ and $\Theta\subset V$ a compatible
lattice. Fix a non-$(\A,\Theta)$-special $\mu\in\vdr$.  Then
there exists an affine arrangement $\B$ with $\B^0=\A^0$ and
representative forms $\fB$, such that each $f\in\catcm$ may be
represented as linear combination of functions from
$\C^\mu_\B[\Tc]$ of the form
\begin{equation}
  \label{eq:trig-parfrac}
   \bet
\prod_{i=1}^n\left(\prod\frac1{\left(1-e_x\right)^{\varepsilon(x)}},\;x\in\fB,
\,x^\circ=\cy_i\right),
\end{equation}
where $t\in\Ts$, $(\cy_1,\dots,\cy_n)\in\bcb(\Ac,\prec)$, and
$\varepsilon:\fB\rightarrow\Z^{\geq0}$ is such that
\[
\sum\{\varepsilon(x)|\; x\in\fB,\,x^\circ=\cy_i\}\neq0,\quad \text{for all }
i\in\fs n.
\] 
\end{prop}

The expression \eqref{eq:trig-parfrac} looks a bit less
satisfactory than \eqref{eq:basis}, because we had to allow
various combinations of products of powers of translates of the
same linear form. There is a further normalization, however,
which allows one to replace \eqref{eq:trig-parfrac} by an element
of $\C^\mu_\B[\Tc]$ of the form
\[ \bet\prod_{i=1}^n \frac1{\left(1-e_{y_i}\right)^{\varepsilon_i}}, \]
where $t\in\Ts$, $(\cy_1,\dots,\cy_n)\in\bcb(\Ac,\prec)$ and
$\varepsilon_i>0$ for $i=1,\dots,n$. We postpone the proof of
this statement to a later publication, as this form of the
expansion is not necessary for our applications here.

Returning to the proof of our theorem, just as in the rational case, a
much weaker statement suffices:
\begin{lem}
   \label{thm:covering} For compatible $\A,\Theta$ with
   $|\Ac|>\dim V=n$ and a nonspecial $\mu\in\vdr$, there is an
   arrangement $\B$ with $\B^\circ=\Ac$, with the property that
   any $f\in\catcm$ may be decomposed into a linear combination
   of elements of $\C^\mu_\B[\Tc]$, where for each element $g$ of
   the decomposition, there is a hyperplane $H_g\in\Ac$, such
   that $g$ is generically regular along {\em any} hyperplane
   parallel to $H_g$.  \end{lem}

\textbf{Step 4.}  We prove the theorem by induction on $|\Ac|$. For
$\catcm$ to be nonempty we need $|\Ac|\geq n$, and the case $|\Ac|=n$
was treated in {\bf Step2}. Thus we will assume $|\Ac|> n$.

It will be useful to formalize the argument at the end of part
{\bf S4} of \thmref{thm:main-rational} in the trigonometric context.To
make our presentation more transparent, we introduce a simplified
notation for the data in
\thmref{thm:main-trig}. The Theorem states the equality of two
quantities, $Z$ and $\widetilde{\ct}$, associated to the data of
$(V,\A,\Gamma,\Theta,f,\tau,\mu)$. Now fix a hyperplane $H\in\A$,
and recall the notation $\amh$ and $\arh$ introduced at the
beginning of {\bf S4}. Assume that $g$ is generically regular
along $H$. Then we have the following 3 sets of data:
\begin{eqnarray*}
d=(V,\A,\Gamma,\Theta,g,\tau,\mu),\\
d\setminus H=(V,\amh,\Gamma,\Theta,g,\tau,\mu),\\
d|H=(H,\arh,\Gamma\cap H,\Theta\cap H,g_{|H},\tau_{|H},\mu_{|H^\circ}).
\end{eqnarray*}
This requires some explanation. First, because of the
compatibility of $H$ and $\Theta$, there are two cases:
\begin{itemize}
\item $H\cap\Theta=\emptyset$;
\item $H\cap\Theta$ is a lattice of full rank in $H$.
\end{itemize}
In the first case the object $d|H$ does not quite make sense, and we
will treat it separately. The second issue is that $H$ is an
affine hyperplane, thus it does not have a canonical vector space
structure. This is a minor technical problem as all of our
constructions are translation invariant. We leave it to the
reader to check that the choice of the origin on $H$ is
immaterial. 

We are ready to state the trigonometric version of
the contraction-deletion principle; in the statement we assume the
notation just introduced.
\begin{lem}\label{thm:del-contr}
If the equality $Z=\widetilde{\ct}$ holds for two of the three
sets of data: $d$, $d\setminus H$ and $d|H$, then it also holds
for the third.
\end{lem}
{\em Proof}: Note that we always have the equality
\[ Z(d) = Z(d\setminus H) - Z(d|H),\]
and our task is to show that the same equality holds for
$\widetilde{\ct}$. The key to this is again \eqref{eq:nbc-ind},
the contraction-deletion relation for \nbc-bases. We impose the
same relation on the ordering of $\Ac$ as we did in {\bf S4}. In
particular, the hyperplane $H^\circ$ is last.

We start with the first case, $H\cap\Theta=\emptyset$, discussed
above. Obviously, we have $Z(d) = Z(d\setminus H)$ here. Analyzing
\eqref{eq:deftodmod}, it is easy to see that the function
$\td_\mu(\Gamma,\ba,\tau)$ vanishes at the common zero $p$ of the set
of affine forms $\ba$ unless $p\in\Gamma$. Clearly, the vertices
$\vert(\At)$ which lie on $H$ are not in $\Gamma$. Then, since $g$ is
regular along $H$ and $p\notin \Gamma$, the iterated constant terms
coming from the image of $s_H$ in \eqref{eq:nbc-ind} all vanish
because of the vanishing of the respective Todd functions at the
vertices. This shows that
\[
\widetilde{\ct}(d) = \widetilde{\ct}(d\setminus H).
\]

In the second case, when $H\cap\Theta$ is a lattice of full rank
in $H$, the proof is analogous to that  of the inductive
statement in part \textbf{S4} of the proof of
\thmref{thm:main-rational}. We will not repeat the reasoning here. The 
only additional condition to check is that
$\mu_{|H^\circ}\in\Delta_{g_{|H}}$. This immediately follows from
the definition \eqref{eq:def-delta}. \qed

Now we are ready to complete the proof of the inductive step. We
start with data $d=(V,\A,\Gamma,\Theta,f,\tau,\mu)$ such that
$|\ac|=N>n$, and assume the equality $Z=\widetilde{\ct}$ for all
cases when $|\Ac|=M<N$. 

We start with applying
\lemref{thm:covering} to our situation. Note that there is a
subtlety here: the condition imposed on $\mu$ in the Theorem  that
$t-\mu$ is not $\Gamma$-special is less restrictive than the
condition in the Lemma, i.e. that $\mu$ is not $\Theta$-special.
This does not cause any problems as both conditions are open,
and we may perturb $\mu$ if necessary. What this means is that
there might be several, essentially different partial fraction
decompositions, which lead to the proof of the same formula
$Z=\widetilde{\ct}$ for a particular data.

Continuing with the proof, denote again the arrangement
guaranteed by the Lemma by $\B$, and let $g$ be one of the terms
of the decomposition. Then, according to
\lemref{thm:covering}, we have $g\in\C^\mu_\mC[\Tc]$ for some
affine arrangement $\mC\subset\B$ with $|\mC^\circ|<N$.  As
$g\in\C^\mu_\mC[\Tc]$, according to the inductive hypothesis
\eqref{eq:main-trig} holds for the data
$(V,\mC,\Gamma,\Theta,g,\tau,\mu)$. Since $g$ is genericly regular 
along all the hyperplanes in $\B\setminus\mC$, we are in
position to use \lemref{thm:del-contr}. This, combined with the
inductive hypothesis, allows us to conclude that the equality
$Z=\widetilde{\ct}$ holds for the data
$(V,\B,\Gamma,\Theta,g,\tau,\mu)$. As this is true for all of the 
terms $g$ in the decomposition of $f$, by additivity, we have
$Z=\widetilde{\ct}$ for
$(V,\B,\Gamma,\Theta,f,\tau,\mu)$. Finally, since $f$ is
genericly regular along the hyperplanes in $\B\setminus\A$, we can use
\lemref{thm:del-contr} again  to conclude the inductive
hypothesis for $(V,\A,\Gamma,\Theta,f,\tau,\mu)$. This completes
the proof.\qed

\subsection{Other forms of the Main Theorem}

Here we rewrite \thmref{thm:main-trig} in a different form. We will
relate the Bernoulli polynomials to the rational trigonometric sums.

First we return to the multiple Bernoulli polynomials introduced in
\secref{sec:polynomial}.
It is clear from \eqref{eq:var-main-rational} that if we represent $f$
as a sum of homogeneous terms, then terms of positive 
degree do not contribute to the Bernoulli polynomial 
$P_{[u]f}^{\A\Gamma}$. This allows us to extend the definition of the 
Bernoulli polynomial to the case of an arbitrary $f\in M_{\A}$.

Using this formalism and comparing Theorems
\ref{thm:var-main-rational} and \ref{thm:main-trig}, we arrive at the
following equality:
\begin{equation} \label{eq:new-form}
Z_f^{\A\Gamma/\Theta}(\tau)=\sum_{p\in \vert(\A/\Theta)}
P_{[u]f}^{\Atp\Gamma}(t), 
\end{equation}
where $t-u\in\Delta_f$. 
\begin{rem} \label{varyu}
This equality may be used to derive various polynomial 
relations among Bernoulli numbers. For example, the left hand 
hand side might vanish for an appropriate $\A$ if 
$|\Gamma/\Theta|$ is small. 

One may simply vary $u$ inside $t-\Delta_f$. When $u$ crosses a wall
of special elements, the terms on the right hand side might change,
but their sum will remain the same.

Another way of finding such relations is to use the freedom in
choosing $t$ such that $e_t=\tau$ when restricted to $\Gamma$. Assume
for simplicity that the toric arrangement $\A/\Theta$ has a single
vertex at 0. Then if $\mu,\mu-s\in\Delta_f$, where $\mu=t-u$ and
$s\in\Theta^*$, then we obtain two essentially different expressions
for $Z_f^{\A\Gamma/\Theta}(\tau)$:
\[ P_{[u]f}^{\A\Gamma}(t) \quad\text{and}\quad
P_{[u]f}^{\A\Gamma}(t+s).\] More generally, we can say that this polynomial will
be constant on the set 
\[\{t+rs|\,t+rs-u\in\Delta_f,\,r\in\Z\}.\]
This is a somewhat surprising "periodicity property" of a linear
combination of Bernoulli polynomials\footnote{The author is grateful
  to Mich\`ele Vergne for pointing out these applications.}.
\end{rem}

Formula \eqref{eq:new-form} reduces the computation of rational
trigonometric sums to that of Bernoulli polynomials, albeit ones
corresponding to meromorphic functions. A further computational
simplification is to represent $f$ locally, near each vertex $p\in
\vert(\At)$ as a product $\mathrm{Pr}_p(f)R_p(f)$, where
$\mathrm{Pr}_p( f)$ is of the form $\prod_{x\in\widehat{\A}_p}
x^{\varepsilon(x)}$ with $\varepsilon(x)\in\Z^{\leq0}$, and $R_p(f)$
is a holomorphic function near $p$, which does not vanish at $p$.
These conditions define $\PP(f)$, which we will call the {\em leading
  principal part}, up to a constant only, but this will be sufficient
for our purposes.

To each linear functional $t$ on $V$ one can associate a first
order linear differential operator $\nabla_t$ on $V^*$: the
directional derivative. Extending this correspondence
multiplicatively, we can associate to every power series $h$,
defined near zero on $V$, a formal differential operator $D^h$ of
possibly infinite order. There is an obvious affine
generalization of this, when one considers affine linear forms
vanishing at a point $p$ and power series near $p$ with a
corresponding operator $D_p^h$. Then, for a nonspecial $t$, we
clearly have
\[ B_f^{\A\Gamma}(t) =   
\left[D^{R(f)}\cdot B_{\mathrm{Pr}(f)}^{\A\Gamma}\right](t),\] 
where the notation means that the differential operator acts on
the function $B$ and then the resulting function is evaluated at
$t$. In this notation we can also write
\[ P_{[u]f}^{\A\Gamma}(t) =
\left[e^{\nabla_\mu}\cdot B^{\A\Gamma}_{f}\right](u),\]
where $u=t-\mu$.

Now it is easy to see that  the equality \eqref{eq:new-form} may be
rewritten in the following form.
\begin{prop}\label{thm:bisform}
  Let $\A,\Gamma\supset\Theta$ be compatible, $u=t-\mu\in V^*$
  nonspecial and $f\in\catcm$. Then we have
\begin{multline}\label{eq:diff-main}
Z^{\A\Gamma/\Theta}_f(\tau)=\sum_{p\in \vert(\A/\Theta)}
\left[D^{R_p(f)}\cdot P^{\Atp\Gamma}_{[u]\PP_{p}(f)}\right](t)=\\
\sum_{p\in \vert(\A/\Theta)} \left[D^{R_p(f)}e^{\nabla_\mu}\cdot
B^{\Atp\Gamma}_{\PP_{p}(f)}\right](t-\mu).
\end{multline}
\end{prop}

This might seem like a clumsy way of presenting our formula, but, 
as it turns out, this is exactly what we need for our application 
in the next section.

\section{The Verlinde formula and the work of Bismut and Labourie}

\subsection{Preliminaries and Verlinde's formula} In this section we
apply our results to a special case which, in fact, motivated our
work.  The computation is related to a natural example of the setup of
the previous section, one provided by Lie theory.   Let $\lig$
be the Lie algebra of $G$. Then in terms of our earlier notation
\begin{itemize}
\item  $V$ is the  Cartan subalgebra of $\lig$;
\item $\ar$ is the HPA on $V$ induced by the set of roots; 
denote by $\har$ the set of positive roots;
\item $\Theta\subset V$ is the coroot lattice;
\item $\Gamma\subset V$ is the lattice for which $\Gamma^*$ is the
  lattice generated by the long roots; for a positive integer $k$ let
  $\Gamma[k]=\{v\in V|\, kv\in\Gamma\}$.
\end{itemize}
There is an inner product $(,)$ on $V^*$ called
\emph{basic}, such that the long roots have square length 2. The
resulting identification of $V$ and $V^*$ induces a one-to-one
map between $\Theta$ and $\Gamma^*$. We will need some additional
notation from Lie theory:
\begin{itemize}
\item Let $\disp\rho=\frac12\sum_{\alpha\in \har}\alpha$,
\item denote by $\theta$ the highest root; the integer
  $h=(\theta,\rho)+1$ is called the dual Coxeter number,
\item the two important functions of Lie theory
\[ d = \prod_{\alpha\in \har} \ala \alpha\quad \text{and}\quad
\delta=\prod_{\alpha\in \har} 2\sqrt{-1}\sin(\pi \alpha), \]
the generalized Vandermonde determinant and the Weyl denominator,
are in $\rfr$ and $\C_\ar[\Tc]$, correspondingly. Note that
$d^{-1}$ is the leading principal part of $\delta^{-1}$ at the origin.
\end{itemize}

It is easy to identify the vertices of our arrangement:
\[ \vert(\ar/\Theta) = \{p+\Theta|\, \langle\alpha\in\har|\,
\alpha(p)\in\Z\rangle_\mathrm{lin}=V^*\}\quad\mathrm{and} \]
\[ \ar_p =\{\{v|\,\alpha(v)=\alpha(p)\}|\,\alpha\in \har,\,\alpha(p)\in\Z\},\]
where $\langle\rangle_\mathrm{lin}$ is the linear span of a set.
Visually, the vertices correspond to the intersections of $n$ or
more root hyperplanes in the Stiefel diagram of the group. For
example, in the case of the rank-2 group $G_2$ there are 5
vertices.

Given positive integers $k,g$ and a dominant weight $\lambda$ such
that $(\theta,\lambda)\leq k$, there is a remarkable formula
discovered by E.~Verlinde \cite{Ver} for a nonnegative integer
$\ver_g(\lambda;k)$, which stands for the dimension of a certain
vector space, the space of ``conformal blocks'' in an appropriate
conformal field theory. It takes the form
\begin{equation}
  \label{eq:verdef}
  \ver_g(\lambda;k) = \frac{((-1)^{|\ar|}|\Gamma/\Theta|(k+h)^n)^{g-1}}{|W|}
\sum_{\gamma\in\Gamma[k+h]/\Theta} \chi_\lambda(\gamma)\delta(\gamma)^{2-2g},
\end{equation}
where $\chi_\lambda$ is the character of the irreducible
representation with highest weight $\lambda$ lifted to the Cartan
subalgebra by the exponential map, and $W$ is the Weyl group of $G$.
The number $\ver_g(\lambda;k)$ vanishes if $\lambda$ is not an integer
linear combination of roots.

Before we proceed, we introduce the shorthand
\begin{itemize}
\item $Z_m^{\ar[k]}(\tau)$ for the rational
trigonometric sum corresponding to the function $\delta^{-m}$,
the lattice $\Gamma[k]$ and some exponential weight $\tau\in\Tc$,
\item $d_p$ for $\prod_{\alpha\in\ar^{+\Theta}_p}(\alpha-\alpha(p))$,
  which is the inverse of the principal part $\PP_p(\delta^{-1})$ at a
  vertex $p$,
\item  $B_m^{\ar_p[k]}$ for the rational sum corresponding
to the arrangement $\ar_p^{+\Theta}$, the lattice $\Gamma[k]$ and
the rational function $d_p^{-m}$, and
\end{itemize}

Recall now the Weyl character formula,
\[\chi_\lambda=\frac1\delta\sum_{w\in W}
\sign(w)e_{w(\lambda+\rho)},\] and the fact that the function
$\delta$ is Weyl-antisymmetric.
Using our notation, we may
write:
\begin{equation}
  \label{eq:ver-v}
\ver_g(\lambda;k) = ((-1)^{|\ar|}|\Gamma/\Theta|(k+h)^n)^{g-1}
Z_{2g-1}^{\ar[k+h]}(\lambda+\rho).
\end{equation}

Now we can apply \thmref{thm:main-trig} to compute the value of 
$\ver_g(\lambda;k)$:
\begin{equation}
  \label{eq:vercomp}
\ver_g(\lambda;k) = (-1)^n((-1)^{|\ar|}|\Gamma/\Theta|(k+h)^n)^{g-1}
\sum_{p\in\vert(\ar/\Theta)} 
\widetilde\ct^{\ar_p\Gamma[k+h]/\Theta}_{\mu,\,\lambda+\rho}
(\delta^{1-2g}),
\end{equation}
where $\mu\in\Delta_{\delta^{1-2g}}$. Clearly,
$\Delta_{\delta^{1-2g}}=(2g-1)\Delta_{\delta^{-1}}$. Denote the
convex polytope $\Delta_{\delta^{-1}}$ simply by $\Delta$ from
now on. The following statement easily follows from the definition 
of the function $\delta$ and \propref{thm:dv}:
\begin{lem} \label{thm:delta}
The set $\Delta$ is the convex hull of the Weyl orbit of
the weight $\rho$. 
\end{lem}

For the brave souls who want to use \eqref{eq:vercomp} for actual
computations, here is the ``reverse engineered'' form:
\begin{multline}
  \label{eq:verdetcomp}
\ver_g(\lambda;k-h) = (-1)^n((-1)^{|\ar|}|
\Gamma/\Theta|k^n)^{g-1}\times\\
\sum_{p\in\vert(\ar/\Theta)} 
\sum_{\ba\in\bB_p} \ict\ba 
  \frac{\sum \{e_{k\tilde t}|\;  
  \tilde t-\lambda/k\in\Gamma^*,\,\tilde t-\mu/k\in\Box(\widehat
  \ba^\circ)\}}{\mathrm{vol}_\Gamma(\widehat \ba^\circ)} 
  \prod_{i=1}^n \frac{k\ala{x}_{i,\widehat \ba}}{{e}^k_{i,\widehat
      \ba^\circ}-1}\frac1{\delta^{2g-1}},
\end{multline}  
where $\vert(\ar/\Theta)$ is defined above and $\bB_p$ is a set of ordered
$n$-tuples of root hyperplanes from $\ar_p$. There are many
choices for $\bB_p$; one possibility is given by 
\eqref{eq:def-nbc}, which depends on a linear ordering of the roots.
We replaced $k$ by $k-h$ to make the formula more readable.

Let us point out once more the computational advantage of the formula
\eqref{eq:vercomp} over the finite sum in \eqref{eq:verdef}: in the
localized formula the number of terms does not depend on $k$, while in
Verlinde's expression the number of terms increases polynomially in
$k$.  The toric arrangement in Example \ref{ex:b2} corresponds to the
case of the group $\mathrm{Spin(5)}$; the lattice $\Gamma$ in the
Example is slightly different from the one appearing in the Lie data,
but this difference causes only minor modifications in the calculations.

\subsection{The Riemann-Roch numbers of the moduli spaces of flat
connections}\label{sec:riemann-roch} 
We will try to give an ultrashort introduction to this subject
here. The reader is referred to \cite{BL,sorger} for more
details. 

We start with the same data as was necessary to describe
Verlinde's formula: we need a simple, simply-connected compact
Lie group $G$, an integer $g$ and an integer $k$ called the {\em
level} and a dominant weight $\lambda$ in the root lattice such
that $(\lambda,\theta)\leq k$. Given these, one may construct a
possibly singular symplectic manifold $\M$, a moduli space of
flat connections on a Riemann surface of genus $g$. This manifold
is endowed with a prequantum line bundle $\LL$. The manifold and
the line bundle depend on the data but we will omit this
dependence from the notation.  Modulo some technical
difficulties, one can define an integer $\chi(\M,\LL)$, the {\em
Riemann-Roch number} of this pair. We will denote this number by
$\chi(g,k,\lambda)$ when we want to emphasize the dependence on
our parameters.

Assume that  holomorphic structures compatible with the
symplectic form are fixed on $\M$ and $\LL$. Then one may define
$\chi(\M,\LL)$ as the alternating sum of dimensions of the sheaf
cohomology groups of the space of holomorphic
sections of $\LL$:
\[ \chi(\M,\LL) =\sum_{i=0}^{\dim \M} \dim H^i(\M,\LL).\]

It is expected that a vanishing theorem holds in this case, which
means that $\dim H^i(\M,\LL)=0$ if $i>0$. On the other hand, the
space of conformal blocks, whose dimension is computed by
Verlinde's formula, may be identified with the space of sections
$H^0(\M,\LL)$ (cf. \cite{blasz}).  Thus one would conjecture that 
\begin{equation}\label{eq:chi=ver}
\chi(g,k,\lambda) = \ver_g(\lambda;k)
\end{equation}

Bismut and Labourie give an explicit formula for the Riemann-Roch
number $\chi(g,k,\lambda)$, and they prove the identity
\eqref{eq:chi=ver} in various cases, notably for large values of $k$
when $\lambda/k$ is fixed. Our goal in this section is to show how our
formalism fits with their formula and to prove \eqref{eq:chi=ver} in
general. Our efforts here are a significant improvement on Section 7
of \cite{BL}.

Note that so far we have described a restricted case of the whole
story, when the ``number of punctures'', denoted by $s$ in
\cite{BL}, is equal to 1. Our proof easily implies the case of
$s>1$ as well {\em if} $g>0$. However, when $s>2$, then $g$, the
genus of the Riemann surface, is allowed to be 0, and we were
{\em not} able to cover this case.

We would like to end this section with a sketchy and incomplete review
of the existing results regarding \eqref{eq:chi=ver}.
\begin{itemize}
\item The equality $\dim H^0(\M,\LL)=\ver_g(\lambda;k)$ is known in
most cases (cf. \cite{blasz,sorger}).
\item The vanishing theorem mentioned above is easy to show in some
cases, and apparently  follows from a recent result of Teleman
\cite{tel}. 
\item The first two points thus provide, albeit a rather rocky road 
to the proof of \eqref{eq:chi=ver}. In a recent paper \cite{MW},
Meinrenken and Woodward prove this equality in complete generality,
including the genus 0 case, using completely different methods.
\end{itemize}

\subsection{The formula of Bismut and Labourie}
Now we are ready to compare our residue theorem to the formula of
Bismut and Labourie for $\chi(g,k,\lambda)$.  We will assume that
the reader has \cite{BL} available for the comparison, since even
introducing all the notation from this reference would have
unreasonably lengthened our paper.

The formula for $\chi(g,k,\lambda)$ is given in Theorems 6.16 and
6.26 of \cite{BL}.  We start with listing the correspondence of
the relevant notation (\cite{BL} $\mapsto$ this paper) in the two
papers:
\begin{itemize}
\item  $(p\mapsto k),\,(c\mapsto h),\,(u\mapsto p),\,(l\mapsto|\ar|)$
\item $(p+c)^{(g-1)\dim(\gz(u))+\frac  s2\dim(\gz(u)/\lt)}\mapsto
  (k+h)^{(2g-2+s)|\ar_p|+n(g-1)}$
\item $\overline{CR}\mapsto\Theta$, $R_l\mapsto\Gamma$,
\item According to \cite[Proposition 1.3]{BL}
  $\mathrm{Vol}(T)^2=|\Gamma/\Theta|$.
\item The function $e_t(p)P_{p,m,q}(t)$ defined by formula (2.166) of
\cite{BL} is equal to $q^{m|\ar_p|}B^{\ar_p[q]}_m(qt)$. This
follows from \eqref{eq:finite}.
\end{itemize}

Consider the case $s=1$ first.  In this case, the formula (6.112)
of \cite{BL} simplifies a little: one may forego the summation over
the Weyl group because of the symmetry.
Then the formula of Bismut and Labourie may be rewritten in
our notation as
\begin{multline} \label{eq:BL}
\chi(g,k,\lambda) =
((-1)^{|\ar|}|\Gamma/\Theta|(k+h)^n)^{g-1} \\
\sum_{p\in
  \vert(\ar/\Theta)} \left[D^{(\delta/d_p)^{1-2g}}e^{\nabla_\mu}
\cdot B_{2g-1}^{\ar_p[k+h]}\right](\lambda+\rho-\mu),
\end{multline}
where $\disp \mu=\rho-\frac{h}{k}\lambda$.

This formula has exactly the same form as the one appearing in
\propref{thm:bisform}.  Taking into account
\eqref{eq:ver-v}, we may conclude
\begin{prop} \label{thm:collect} 
Let $k,g>0$, $s=1$. Then 
$ \chi(g,k,\lambda) = \ver_g(\lambda;k)$ as long as $\disp
\mu=\rho-\frac{h}{k}\lambda\in \Delta$ and
$\disp \frac\lambda k$ is not $\Gamma$-special.
\end{prop}
Note that here we wrote $\Delta$ instead of $(2g-1)\Delta$, i.e. we
set $g=1$, as this case clearly implies the cases of higher genera.

Geometricly, the condition that $\lambda/k$ is special means that the
moduli space $\M$ is singular. Let us write \eqref{eq:BL} in the form
\eqref{eq:new-form}. We obtain
\[ \chi(g,k,\lambda) =
((-1)^{|\ar|}|\Gamma/\Theta|(k+h)^n)^{g-1} 
\sum_{p\in
  \vert(\ar/\Theta)} P^{\ar_p\Gamma[k+h]}_{[u]
\delta^{1-2g}}(\lambda+\rho),
\] %
where $\disp u=(k+h)\frac{\lambda}k$. Now if $\disp \frac\lambda
k$ is special, then we can use \remref{varyu} and conclude that
by perturbing $u$ here a bit, we still have a valid formula. This
perturbation corresponds to computing the Riemann-Roch number on
a different moduli space, which is smooth. This moduli space may
be different if $u$ moves in different directions, but
\remref{varyu} tells us that we will always get the same answer
(cf. \cite[Theorem 6.40]{BL}). A similar trick may be used if
$\mu$ is on the boundary of $\Delta$. The reader is referred to
\cite[Section 6.11]{BL} for further explanations.

The case of multiple punctures, $s>1$, easily reduces to the same
argument. Indeed, each term in the sum over the $s$ copies of the
Weyl group in formula (6.112) of \cite{BL} is of the form of a
product of factors of the kind that we encounter in the $s=1$,
$g>0$ case. Thus if the condition in the Proposition holds, then we
have $\mu_i\in \Delta_{\delta^{-1}}$ for each $i=1,\dots,s$ and
then we can use the easy direction of \propref{thm:dv} to conclude
that the appropriate condition $\mu\in\Delta$ holds in this case.

When $s>2$ and $g=0$, then the number of punctures exceeds the
number of factors of $\delta$ in the denominator, and this argument
would require $\mu\in\frac13\Delta$, which does not always hold.

\subsection{Studying the condition} Now we turn to the study of
the condition in \propref{thm:collect}. We maintain the notation
of the previous section. Recall that $\lambda$ is a dominant
weight satisfying $(\theta,\lambda)\leq k$. We denoted the
polytope $\Delta_{\delta^{-1}}$ simply by $\Delta$, and according
to \lemref{thm:delta}, $\Delta$ is the convex hull of the
$W$-orbit of $\rho$.  Then the condition in \propref{thm:collect}
may be rewritten as follows.
\begin{prop}
  \label{thm:tech}
Let $\ar$ be a root system of a simple Lie algebra with a chosen
dominant chamber $\dc\subset V^*$. Then if $(\theta,\nu)< h$ and
$\nu\in\dc$, then $\rho-\nu\in\Delta$.
\end{prop}
\emph{Proof}: The statement can be proved by direct computation using
the classification of simple Lie algebras. We will demonstrate the
method in the cases of the Lie algebras $A_n$ and $D_n$. Note that the
lattices associated to the Lie algebra do not figure in this
statement; we only need to know the set of roots to formulate it.

Denote the simple roots by $\{\alpha_i\}_{i=1}^n$; thus
$\dc=\{\nu|\,(\alpha_i,\nu)\geq0,\,i\in\fs n\}$. Denote by
$\omega_i$, $i\in\fs n$, the corresponding fundamental weights
rescaled as follows: $(\alpha_i,\omega_j)=0$ for $i\neq j$, and
$(\theta,\omega_i)=h$.  Then the set
$\dch=\{\nu\in\dc|\,(\theta,\nu)\leq h)\}$ is a simplex with
vertices at $0,\omega_1,\dots\omega_n$. As $\Delta$ is the
convex hull of the $W$-orbit of $\rho$, it is clear that to prove the
statement, one needs to show that $\rho-\omega_i\in\Delta$ for
$i\in\fs n$.

For the Lie algebra $A_n$ this is true in a remarkable
fashion. Here we may think of $V^*$ as the subspace in the
Euclidean $n+1$-space with coordinate vectors
$\delta_i,\,i=0,\dots,n$. To simplify the notation, let
$v=\sum_{i=0}^{n}\lambda_i\delta_i=(\lambda_0,\lambda_1,\dots,\lambda_n)$
be a generic vector in this space. Then
\begin{itemize}
\item $V^*=\{v|\,\lambda_0+\lambda_1+\dots+\lambda_n=0\}$,
\item $\dc=\{v\in V^*|\,\lambda_0\geq \lambda_2\geq\dots\geq\lambda_n\}$,
\item $\rho = \frac12\smn (n-2i)\delta_i$,
\item $\theta=\delta_0-\delta_n$, $h=n+1$ and
\item
  $\omega_m=m\sum_{i=0}^{n-m}\delta_i-(n+1-m)\sum_{i=n-m+1}^n\delta_i$,
  $i=1,\dots,n$.
\end{itemize}
Computing the coordinates of the vector $\rho-\omega_m$ for each
$m$, one discovers that they are the same as those of $\rho$, but
in different order. Since the Weyl group of $A_n$ acts by
permuting the coordinates this means that the vertices of the
simplex $\rho-\dch$ are a {\em subset of the vertices} of the
convex polytope $\Delta$. This is, of course, a much stronger
statement than what we needed, but it also shows that our
computation is ``sharp'' in a certain sense.

This sharp statement does not hold for the other simple Lie
algebras. The general method of checking the condition goes as
follows. For $v\in V^*$, denote by $v^\#$ the vector in $\dc$ which
is Weyl equivalent to $v$. Since the codimension-1 faces of the
polytope $\Delta$ lie in hyperplanes perpendicular to Weyl shifted
fundamental weights, and $\Delta$ is Weyl invariant, we have
\[
\Delta=\{v\in
V^*|\,(v^\#,\omega_i)\leq(\omega_i,\rho),
\,i\in\fs{n}\}
\]
Hence in order to prove the Proposition for $\ar$, one needs to
check the $n^2$ inequalities:
\begin{equation}\label{eq:nn}
((\rho-\omega_m)^\#,\omega_i)\leq(\omega_i,\rho),\quad i,m\in\fs n.
\end{equation}
Here is what happens in the case of $D_n$. We have
\begin{itemize}
\item $V^*=\{v=\sum_{i=1}^n\lambda_i\delta_i\}$ is the Euclidean $n$-space,
\item $\dc=\{v\in V^*|\,\lambda_1\geq
\lambda_2\geq\dots\geq|\lambda_n|\}$, 
\item The $j$th coordinate of $v^\#$ is the $j$th largest
number among the absolute values of the coordinates of $v$, with
the possible exception of the $n$th coordinate, which has a minus
sign if the number of negative coordinates of $v$ is odd.
\item $\rho = \sum_{i=1}^n (n-i)\delta_i$,
\item $\theta=\delta_1+\delta_2$, $h=2n-2$,
\item $\omega_1=2(n-1)\delta_1$, $\omega_m=(n-1)\sum_{i=1}^m\delta_i$
  for $m\in\overline{2,n-2}$,
  $\omega_{n-1}=\omega_{n-2}+(n-1)(\delta_{n-1}+\delta_n)$, and
$\omega_{n}=\omega_{n-2}+(n-1)(\delta_{n-1}-\delta_n).$
\end{itemize}
The inequalities \eqref{eq:nn} are easy to check now. Indeed, for
example, for $m\in\overline{2,n-2}$ we have 
\[ \rho-\omega_m = \sum_{i=1}^m (1-i)\delta_i+\sum_{i=m+1}^n
(n-i)\delta_i. \] 
In this case, we see that comparing the $i$th
coordinates, we have 
\[ |(\delta_i,\rho)|\geq|(\delta_i,\rho-\omega_m)|,\] which
immediately implies the inequalities. 

Note that while for $A_n$ we have $(\rho-\omega_m)^\#=\rho$ for all
$m\in\fs n$, here this only holds for $m=1,n-1,n$. The cases of the
Lie algebras $B_n$ and $C_n$ are completely analogous and are left as
an exercise to the reader. The exceptional Lie algebras need to be
checked on a case by case basis. We will not present these tedious
calculations here; rather, we hope that someone will find a
conceptual Lie theoretic proof of the Proposition, which will not rely
on the classification theorem. \qed

Let us summarize our results. If we assume that
$\chi(g,k,\lambda)$ is the the appropriately modified definition
of the Riemann-Roch number of the moduli spaces (cf. the
discussion after \propref{thm:collect}), then \propref{thm:tech}
combined with \propref{thm:collect} implies
\begin{thm}\label{thm:chi=ver} 
 For arbitrary $k$ and $\lambda$,
 we have
\[ \chi(g,k,\lambda)=\ver_g(\lambda;k) \]
whenever $g>0$.
\end{thm}

As we pointed out earlier, our arguments prove the same statement
in the case of an arbitrary number of punctures as well. The reader is
referred to \cite{MW} for a different approach to this Theorem which also
covers the genus 0 case.

\subsection{Quasipolynomial behavior and the topology of moduli
spaces} As we mentioned in the introduction, one of the reasons
for searching for a localization formula for Verlinde's function
$\ver_g(\lambda,k)$ was to uncover the polynomial nature of this
function. The coefficients of this polynomial, in turn, give information about the topology of
the moduli spaces in term of certain characteristic numbers. This
point, first raised in \cite{thad}, was elaborated upon in 
\cite{comb}, mainly with the group $SU(n)$ in mind.  As explained
in the paper of Bismut and Labourie \cite{BL}, the situation is
more complicated for other groups: the Riemann-Roch Theorem needs
to be replaced by the Kawasaki-Riemann-Roch Theorem for
orbifolds.

As we will see, in this case instead of polynomials one deals
with quasipolynomials.  We say that a function on the positive
integers $q(k)$ is a {\em
  quasipolynomial} in $k$ of order $L$ if the function $q(kL+a)$
is a polynomial in $k$ for every $a\in\Z$. Equivalently, a quasipolynomial
of order $L$ is a linear combination of characters of the cyclic
group $\Z/L\Z$ with values in polynomials.

In this concluding paragraph, we would like to explain how one
may obtain information about the components of the
quasipolynomial behavior. This might be
useful in problems such as determining the Picard group of the
moduli space (cf. \cite{picard}). Such computations also serve as
a handy mnemonic for understanding the stratification of the
moduli spaces worked out in \cite{BL}. We would like to emphasize that
the statements detailed below may be read off from the results of
\cite{BL} already. We only add a somewhat more compact and transparent
formalism for the computations.

We start by noting a few simple scaling properties of the Bernoulli
functions:
\begin{itemize}
\item For the case when the vertex $p$ is at the origin, we have
\[ B_m^{\ar[k]}(kt) = k^{m|\ar|}B_m^{\ar[1]}(t).\]
\item If the vertex is not at the origin, but we have
$p\in\Gamma$, then this identity is generalized as follows:
\begin{equation}\label{eq:bk}
 B_m^{\ar_p[k]}(kt) = e_t(kp)k^{m|\ar_p|}B_m^{\ar_p[1]}(t).
\end{equation}
\item If the vertex $p$ is not in $\Gamma$, then the scaling
  properties are a bit more complicated. However, for any vertex
  $p\in\vert(\ar^{+\Theta})$ we have $Np\in\vert(\ar^{+\Theta})$ for
  all integers $N$; also, we have $Mp\in\Gamma$ for some integer $M$.
  In particular, $\Gamma\subset\ar^{+\Theta}$.
\end{itemize}

As a next step, we amend the material of \secref{sec:riemann-roch}
slightly. We need to make the relation between the parameters
$g,k,\lambda$ and the prequantized moduli spaces $(\M,\LL)$ more
precise.
\begin{itemize}
\item  The space $\M(g,k,\lambda)$ depends on $g$ and the {\em quotient}
  $\lambda/k$ only.  Thus we may write $\M(g,\lambda/k)$. 
\item The $m$th power of the line bundle $\LL(g,k,\lambda)$, defined on
  $\M(g,\lambda/k)$, is simply the bundle $\LL(g,mk,m\lambda)$
\item One may, in fact, define $\M(g,t)$ for any $t\in\dc^{\leq1}$;
  these spaces form a smooth family, when $t$ varies inside a chamber
  of nonspecial elements. 
\end{itemize}
For a given set of data $(g,k,\lambda)$, we are interested in
computing the orbifold Riemann-Roch number
\begin{equation}
  \label{eq:quasi}
\chi(g,km,m\lambda)\int_{\M} e^{m\LL}\td(\M)  
\end{equation}
as a function of $m$. Analyzing \eqref{eq:BL}, we see that $\mu=\rho-ht$,
where $t=\lambda/k$, 
and thus does not depend on $m$. Then the formula has the general form
\[ \chi(g,mk,m\lambda) = \mathrm{const}\;(mk+h)^{n(g-1)}\sum_{p\in
  \vert(\ar/\Theta)} D_p\cdot B_{2g-1}^{\ar_p[mk+h]}((mk+h)t),\]
where $D_p$ is an operator that does not depend on $m$. Naturally,
here we assumed that $t$ is not special. If $t$ is special, the we
need to use the formalism of the Bernoulli polynomials $P_{[u]}$.

Comparing this formula for $\chi$ to the scaling properties of the
Bernoulli functions that we listed above, we arrive at the following
general conclusions.
\begin{prop}
  \label{thm:last}
  The vertices $\vert(\ar/\Theta)\subset T$ are a union of cyclic
  subgroups, one of which is $\Gamma/\Theta$. The orbifold
  Riemann-Roch number $\chi(g,mk,m\lambda)$ is a polynomial in $m$ if
  $kp\in\Gamma$ and $\lambda(p)\in\Z$ for all
  $p\in\vert(\ar^{+\Theta})$.
\end{prop}
For example, in the case of the Lie group $G_2$ the set of vertices
$\vert(\ar/\Theta)$ consists of 5 points: the union of a cyclic group
of order 3, which is $\Gamma/\Theta$, and a two cyclic groups of order
2. This means that if we take the line bundle corresponding to the
trivial character ($\lambda=0$) at level 1 ($k=1$), then the orbifold
Euler characteristic will be a quasipolynomial of order 2. If we take
$\lambda\neq0$, then $\chi(g,mk,m\lambda)$ could end up being a
quasipolynomial of order 6.


\newpage
\begin{thebibliography}{99}
\bibitem{AB} M. Atiyah, R. Bott, { \em The Yang-Mills equations over Riemann surfaces}, Philos. Trans. Roy. Soc. London Ser. A 308 (1983), no. 1505, 523--615.


\bibitem{blasz} A. Beauville, Y. Laszlo, {\em Conformal blocks and
generalized theta functions}, Comm. Math. Phys. {\bf 164} (1994),
no. 2, 385--419.

\bibitem{BL} J.-M. Bismut, F. Labourie, {\em Symplectic geometry
and the Verlinde formulas.} Surveys in differential geometry:
differential 
geometry inspired by string theory, pp. 97--311, Surv. Differ. Geom.,
{\bf 5}, Int. Press, Boston, MA, 1999

\bibitem{BVvp} M. Brion, M. Vergne, {\em Residue formulae, vector
partition functions and lattice points in rational
polytopes}, J. Amer. Math. Soc. {\bf 10} (1997), no. 4, 797--833.

\bibitem{bv} M. Brion; M. Vergne, {\em Arrangement of
hyperplanes. I. Rational functions and Jeffrey-Kirwan
residue}, Ann. Sci. École Norm. Sup. (4) {\bf 32} 
(1999), no. 5, 715--741.

\bibitem{BV} M. Brion, M. Vergne, {\em Arrangement of
hyperplanes. II. The Szenes formula and Eisenstein series}, Duke
Math. J. {\bf 103} (2000), no. 2, 279--302.

\bibitem{zel} I.~M.~Gelfand, A.~V.~Zelevinski, {\em Algebraic and
combinatorial aspects of the general theory of hypergeometric
functions.}  (Russian) Funktsional. Anal. i Prilozhen. {\bf 20}
(1986), No.3, 17--34.

\bibitem{JK} L. Jeffrey, F.  Kirwan, {\em Intersection theory on
moduli spaces of holomorphic bundles of arbitrary rank on a
Riemann surface}, Ann. of Math. 
(2) {\bf 148} (1998), no. 1, 109--196.


\bibitem{picard} S. Kumar, M. S. Narasimhan, {\em Picard group of the moduli
spaces of $G$-bundles} Math. Ann. {\bf 308} (1997), no. 1, 155--173.

\bibitem{MW} E. Meinrenken, C. Woodward, {\em Hamiltonian loop
group actions and Verlinde factorization}, J. Differential
Geom. {\bf 50} (1998), no. 3, 417--469.

\bibitem{OT} P. Orlik, H. Terao, {\em Arrangements of hyperplanes}, Springer-Verlag, Berlin, 1992

\bibitem{sorger} C. Sorger, {\em La formule de Verlinde}, Séminaire Bourbaki, Vol. 1994/95. Astérisque No. 237 (1996), Exp. No. 794,
3, 87--114.

\bibitem{comb} A. Szenes, {\em The combinatorics of the Verlinde
formulas}, in {\sl Vector bundles in algebraic geometry (Durham,
1993)}, 241--253, London Math. Soc. Lecture Note Ser.,
208.

\bibitem{bern}
A. Szenes \emph{Iterated Residues and multiple Bernoulli polynomials}
 Internat. Math. Res. Notices 1998, no. 18, 937--956. (\textsl{arXiv:hep-th/9707114}).

\bibitem{SzV} A. Szenes, M. Vergne, in preparation.

\bibitem{tel} C. Teleman, {\em The quantization conjecture
revisited}, Ann. of Math. (2) {\bf 152} (2000), no. 1, 1--43.  

\bibitem{thad} M. Thaddeus, {\em Conformal field theory and the
cohomology of the moduli space of stable bundles}, J. Differential
Geom. {\bf 35} (1992), no. 1, 131--149.

\bibitem{Ver} E. Verlinde, {\em Fusion rules and modular
transformations in $2$D conformal field theory}, Nuclear Phys. B
{\bf 300} (1988), no. 3, 360--376.

\bibitem{wit1} E. Witten, {\em On quantum gauge theories in two
dimensions} Comm. Math. Phys. {\bf 141} (1991), no. 1, 153--209.

\bibitem{wit2} E. Witten, {\em Two-dimensional gauge theories revisited},
J. Geom. Phys. {\bf 9} (1992), no. 4, 303--368.


\end{thebibliography}
\end{document}